\documentclass{amsart}

\usepackage[english]{babel}
\usepackage{graphicx}
\usepackage[hidelinks]{hyperref} 
\usepackage{amssymb}
\usepackage{amsmath}
\usepackage[foot]{amsaddr}
\usepackage{amsfonts}
\usepackage[dvipsnames, table]{xcolor}
\usepackage{mathtools}
\usepackage{algpseudocode}
\usepackage{algorithm}
\usepackage{comment}
\usepackage[top=1in, bottom=1.25in, left=1.25in, right=1.25in]{geometry}
\usepackage{float}
\usepackage{tikz}
\usepackage{multirow}
\usepackage[position=top]{subfig}
\usepackage{epstopdf}
\usepackage[shortlabels]{enumitem}
\newcommand{\first}[1]{{\color{black}{#1}}}
\newcommand{\second}[1]{{\color{black}{#1}}}
\newcommand{\both}[1]{{\color{black}{#1}}}
\newcommand{\norm}[1]{\lVert #1 \rVert}

\newtheorem{example}{Example}{}
\newtheorem{definition}[example]{Definition}
\newtheorem{cor}[example]{Corollary}{}
\newtheorem{thm}[example]{Theorem}{}
\newtheorem{remark}[example]{Remark}
\newtheorem{assump}[example]{Assumption}
\newtheorem{proposition}[example]{Proposition}{}

\numberwithin{equation}{section}
\numberwithin{example}{section}

\graphicspath{{JNMA/}}

\usepackage{todonotes} 
\newcommand{\ms}[1]{{ \color{magenta}{[MS - #1]}}}
\newcommand{\ls}[1]{{ \color{teal}{[LS - #1]}}}

\begin{document}

\title[DLRA strategies for nonlinear feedback control]{Dynamical Low-Rank Approximation Strategies for Nonlinear Feedback Control Problems}

\author[L. Saluzzi and M. Strazzullo]{Luca Saluzzi$^{1,*}$ and Maria Strazzullo$^{2,*}$}
\date{\today}

\address{
$^1$ Dipartimento di Matematica, Sapienza Universit\`a di Roma, Piazzale
Aldo Moro, 5, Roma, 00185, Italy.
}
\address{
$^2$ Politecnico di Torino, Department of Mathematical Sciences ``Giuseppe Luigi Lagrange'', Corso Duca degli Abruzzi, 24, 10129, Torino, Italy.}
\address{
$^*$ INDAM, GNCS member.}
\email{luca.saluzzi@uniroma1.it, maria.strazzullo@polito.it}
\begin{abstract}

This paper addresses the stabilization of dynamical systems in the infinite horizon optimal control setting using nonlinear feedback control based on State-Dependent Riccati Equations (SDREs). While effective, the practical implementation of such feedback strategies is often constrained by the high dimensionality of state spaces and the computational challenges associated with solving SDREs, particularly in parametric scenarios. To mitigate these limitations, we introduce the Dynamical Low-Rank Approximation (DLRA) methodology, which provides an efficient and accurate framework for addressing high-dimensional feedback control problems. DLRA dynamically constructs a compact, low-dimensional representation that evolves with the problem, enabling the simultaneous resolution of multiple parametric instances in real-time.

We propose two novel algorithms to enhance numerical performances: the cascade-Newton-Kleinman method and Riccati-based DLRA (R-DLRA). The cascade-Newton-Kleinman method accelerates convergence by leveraging Riccati solutions from the nearby parameter or time instance, supported by a theoretical sensitivity analysis. R-DLRA integrates Riccati information into the DLRA basis construction to improve the quality of the solution. These approaches are validated through nonlinear one-dimensional and two-dimensional test cases showing transport-like behavior, demonstrating that R-DLRA outperforms standard DLRA and Proper Orthogonal Decomposition-based model order reduction in both speed and accuracy, offering a superior alternative to Full Order Model solutions.

\end{abstract}

\maketitle

\section{Introduction}
\label{sec:intro}
Optimal control of dynamical systems is a central topic in diverse scientific fields. Indeed, it is ubiquitous in many engineering and industrial applications based on the stabilization of nonlinear dynamical systems toward a dynamic equilibrium.
Nonlinear stabilization is commonly addressed through the implementation of feedback (closed-loop) controllers. Unlike open-loop controls, these controllers provide superior stability 
features in the presence of external disturbances.
The design of optimal feedback controls typically relies on dynamic programming, which formulates the optimal feedback policy as the solution to a Hamilton-Jacobi-Bellman (HJB) nonlinear Partial Differential Equation (PDE).
The progress in the solution of HJB PDEs for optimal control is significant. For example, substantial attention has been drawn to sparse grids strategies \cite{GK16}, tree-structure approaches \cite{falcone2023approximation}, applications of artificial neural networks \cite{Darbon_Langlois_Meng_2020,Kunisch_Walter_2021,Zhou_2021} and regression-type methods
in tensor formats \cite{oster2022approximating,dolgov2023data}.

In this study, we address the problem of feedback stabilization by employing an alternative methodology that integrates concepts from dynamic programming and \first{Model Predictive Control (MPC)}. Specifically, we utilize the State-Dependent Riccati Equation (SDRE) approach \cite{ccimen2008state}. Unlike traditional methods that require solving the HJB equation, which is often computationally prohibitive for nonlinear systems, the SDRE approach offers a practical and computationally efficient alternative. By reformulating the control problem into a state-dependent framework, the SDRE method enables the synthesis of near-optimal feedback controllers while avoiding the 
numerical burden associated with the direct solution of the HJB equation. \first{At each sampling instant, the controller relinearizes the dynamics around the current state, solves an algebraic Riccati equation, and applies the resulting gain, thereby emulating MPC’s update-and-apply cycle.
Unlike MPC, however, SDRE does not enforce state or input constraints natively and requires only a single algebraic Riccati solve per step rather than a full finite-horizon nonlinear program.}

For high-dimensional Riccati equations, effective \emph{ad hoc} solvers are available, such as those described in \cite{Kirsten_Simoncini_2020,BBKS_2020}.  In this contribution, we employ the Newton-Kleinman (\texttt{NK}) strategy, an iterative technique that solves the linearized Riccati equation at each step until a specified convergence criterion is achieved \cite{kleinman1968iterative}. Despite the availability of advanced iterative solvers, the application of the SDRE approach remains challenging, particularly in the context of time-dependent parametric settings. These scenarios require solving numerous high-dimensional Riccati equations across a wide range of parameters and time instances, posing significant computational demands for control design and implementation. 
Nevertheless, nonlinear control theory can be of utmost importance in real-time and many query contests, despite the limited applicability due to the burden of the computational costs. For this reason, Reduced Order Models (ROMs) for control problems have recently gained increasing attention in the scientific community. ROMs build a low-dimensional framework spanned by basis functions related to properly chosen parametric instances in an \emph{offline stage}. Once provided the reduced space, a novel evaluation of the parameter is performed by a Galerkin projection in an \emph{online phase}. In this way, ROMs are faster than the Full Order Model (FOM) upon which they are built, and solve the system with a certain degree of accuracy. Here, we report a non-exhaustive list of contributions to the field based on linear projection-based ROM for nonlinear control in open-loop \cite{hinze2005proper,Kunisch1999345,PichiStrazzullo,Strazzullo3,StrazzulloVicini} and closed-loop contexts
\cite{KVX2004,HV2005,alla2020feedback,kirsten2024multilinear,de2023optimal}.
The presented strategies can provide a robust and efficient representation, however, they might fail in transport and wave-like settings, where a large number of basis functions is needed to achieve an acceptable accuracy \cite{chapKolmogorov}. Many strategies have been proposed to overcome this issue, the interested reader may refer to the complete literature overview proposed in \cite{Peher} and the references therein.

In this work, we deal with parametric nonlinear infinite horizon optimal control problems, using a nonlinear projection-based ROM whose basis functions evolve over time and follow the system dynamics. 
We exploit the Dynamical Low-Rank Approximation (DLRA) \cite{KocLubich}, introduced to efficiently compute approximations of large-scale time-dependent problems exploiting time integration in the tangent space of the reduced representation of the solution. DLRA provides a \first{low-rank} representation of the system and the evolution of the basis ensures an increased accuracy if compared to standard data-compression techniques such as Proper Orthogonal Decomposition (POD), see, e.g., \cite{sapsis2009dynamically,Kusch20221} for stochastic PDEs and \cite{Pagliantini2021409} for Hamiltonian systems.
The application of the DLRA methodology in the context of optimal control and the HJB equation was previously explored in \cite{eigel2023dynamical}, where the authors addressed the approximation of finite horizon optimal control problems. Their approach leverages policy iteration and a dynamical low-rank approximation in tensor format for the value function.
In contrast, this work focuses on an infinite horizon parametric setting, solved using the SDRE technique. Here, the DLRA is employed to achieve a \first{low-dimensional} and efficient representation of the solution to the parametric dynamical system.
Moreover, we propose novel strategies to (i) accelerate the solution of SDREs arising from time integration and (ii) enhance the accuracy of DLRA by incorporating information on the control action and the SDRE solution. This leads to the development of a Riccati-based DLRA (R-DLRA), which adapts concepts from \cite{schmidt2015basis, schmidt2018reduced} to further refine the methodology.

The originality of this contribution relies on:
\begin{itemize}
    \item a novel strategy to solve iteratively the SDRE guided by a theoretical sensitivity analysis: the cascade-Newton-Kleinman (\texttt{C-NK}) method. The \texttt{C-NK} exploits the parametric and evolutive information of the system to accelerate the convergence of the classical \texttt{NK}, providing a faster tool to solve Riccati problems both at the FOM and ROM level;
    \item the application of the DLRA approach combined to \texttt{C-NK} for infinite horizon optimal control problems solved by SDRE. The DLRA proves to be a fast and efficient method for representing controlled dynamics in scenarios that pose significant challenges for traditional linear projection-based ROMs;
    \item the enhancement of the classical DLRA approach employing the Riccati information denoted as R-DLRA. The basis functions of the DLRA are enriched with basis functions related to the SDRE or the control action. These additional basis functions increase the accuracy of the DLRA and serve as a tailored strategy for control problems.
\end{itemize}
The DLRA and R-DLRA algorithms provide a fascinating alternative to POD to address the \first{\emph{curse of dimensionality} resulting from spatial discretizations in control applications governed by transport phenomena.}
We validate our findings with one-dimensional and two-dimensional nonlinear \first{PDEs, where the terms "one-dimensional" and "two-dimensional" refer specifically to the spatial dimensionality of the computational domain.}

The paper is outlined as follows: in Section \ref{sec:problem}, we introduce the parametric infinite horizon feedback control problem we address, the SDRE deriving from the control system, its numerical approximation, and the novel \texttt{C-NK} approach. Section \ref{sec:sensitivity} provides a theoretical justification for the employment of the \texttt{C-NK}. Section \ref{sec:DOrb} introduces DLRA and its extension to feedback control. In Section \ref{sec:algorithm} we propose two novel R-DLRA strategies based on the enhancement of DLRA with Riccati information and \first{the \texttt{C-NK}} approach. The numerical algorithms are validated in Section \ref{sec:results}, where we compare the performance of DLRA, R-DLRA, and POD in terms of computational efficiency and accuracy against the FOM solution for feedback control problems based on the Burgers' equation in one and two dimensions. Conclusions follow in Section \ref{sec:conclusions}.

\section{Problem Formulation}
\label{sec:problem}
This section introduces the parametric infinite horizon nonlinear feedback control problem we deal with in this contribution. The problem formulation is described in Section \ref{sec:problem_ihfc}, while the SDRE method is described in Section \ref{sec:SDRE}. Finally, in Section \ref{sec:NK}, we discuss the numerical approximation of SDREs.

\subsection{{Parametric} infinite horizon feedback control}
\label{sec:problem_ihfc}
We {now focus on} 
asymptotic stabilization by infinite horizon optimal control {in a parametric setting}. Given a parameter $\mu \in \mathcal P \subset \mathbb R^{n}$, with $n \in \mathbb N$, the main objective is to minimize the following \emph{cost functional} \begin{equation}
\label{eq:cost_fun}
    \min_{u(\cdot) \in \mathcal{U}_{ad}} J(u(\cdot),y_0; \mu) =  \int_0^\infty y(s; \mu)^\top Q\, y(s; \mu) + u(s; \mu)^\top R \,u(s; \mu) \;{ds},
\end{equation}
subject to a nonlinear dynamics
\begin{equation}
\label{eq:dynamic}
    \dot{y}(t; \mu) = f(y(t; \mu)) + B u(t; \mu), \quad y(0;\mu) = y_0(\mu), 
\end{equation}
where $y(t; \mu)$ is the \emph{state} variable represented by {a} column vector in $\mathbb R^{N_h}$, $u(\cdot; \mu) \in \mathcal{U}_{ad} = \{\second{v}: \mathbb R^+ \rightarrow U \subset \mathbb R^m, \; \text{measurable}\}$ {is the \emph{control} variable}, $f: \mathbb{R}^{N_h} \rightarrow \mathbb{R}^{N_h}$ is the system dynamics, $B \in \mathbb{R}^{N_h \times m}$ is the control matrix, $Q \in \mathbb R^{N_h \times N_h}$ is a symmetric \both{positive} semidefinite matrix, $R \in \mathbb R^{m \times m}$ is a symmetric positive definite matrix, {and $y_0(\mu)$ is a parametric initial condition.}
For the sake of notation, we write the state and the control in a more compact form as $y_{\mu,t} = y(t; \mu)$ and $u_{\mu,t} = u(t; \mu)$. Following the same notation, the parametric initial condition will be $y_0(\mu) = y_{\mu,0}$. We will use both the extended and the compact notations, depending on the needs.

One method for the resolution of the optimal control problem \eqref{eq:cost_fun}-\eqref{eq:dynamic} involves solving the 
HJB equation, a first-order nonlinear PDE defined in $\mathbb{R}^{N_h}$. 

In many applications, the dimension $N_h$ may be large and yields the \emph{curse of dimensionality}, a phenomenon where the computational resources required to solve the HJB grow exponentially with $N_h$.
This 
is a significant limitation of traditional numerical methods. Standard grid-based approaches, for instance, become computationally infeasible as the required number of grid points increases exponentially in higher dimensions. 
In contrast to these computationally intensive approaches, a widely used alternative in practical applications is the SDRE method, which balances computational efficiency and control performance in high-dimensional settings.

\subsection{State-Dependent Riccati Equation}
\label{sec:SDRE}

The SDRE 
is a 
{powerful} mathematical 
{tool} with broad-ranging practical uses \cite{ccimen2008state}. As a natural extension of the traditional Riccati equation, the SDRE expands its versatility by integrating state-dependent coefficients, thereby adapting to nonlinear and fluctuating systems. This method hinges on sequentially solving linear-quadratic control tasks derived from a step-by-step linearization of the dynamics. 
\second{Assume that \( f \in C^1(\mathbb{R}^{N_h}) \) and satisfies \( f(0) = 0 \). Then, according to Proposition~1 in~\cite{cimen2010systematic}, there exists a state-dependent matrix function  $A(x): \mathbb{R}^{N_h} \rightarrow \mathbb{R}^{N_h \times N_h}$ such that $f(x)=A(x)x$. Consequently, the dynamical system \eqref{eq:dynamic} can be rewritten as}

\begin{align*}
    \dot{y}(t; \mu) & = A(y(t; \mu)) y(t; \mu) +B u(t; \mu), \\
\second{y(0;\mu)} &= y_{\mu,0}.
    \label{semilinear}
\end{align*}
Assuming that the system dynamics are linear with respect to the state, i.e., $A(y(t)) = A\in\mathbb{R}^{N_h\times N_h}$, the problem reduces to what is commonly known as the Linear Quadratic Regulator (LQR) problem. In this context, we introduce the notions of stabilizability and detectability \first{for some general matrices $A \in \mathbb{R}^{N_h\times N_h}, B \in \mathbb{R}^{N_h\times m}$ and $ C\in \mathbb{R}^{N_h\times N_h}$}, which will play a crucial role in the subsequent analysis.

\begin{definition}
   \first{Let us consider $A \in \mathbb{R}^{N_h\times N_h}$ and $B \in \mathbb{R}^{N_h\times m}$}. The pair $(A,B)$ is stabilizable if there exists
a feedback matrix \first{$K \in \mathbb{R}^{m \times N_h}$} such that $A -B K$
is stable, i.e., if all its eigenvalues lie on the
open left half complex plane. 
\end{definition}


\begin{definition}
    \first{Let us consider $A \in \mathbb{R}^{N_h\times N_h}$ and $C \in \mathbb{R}^{N_h\times N_h}$}. The pair $(A,C)$ is detectable if the pair $(A^\top, C^\top)$ is stabilizable.

\end{definition}

\first{In the LQR setting,} we suppose that the pair $(A,B)$ is stabilizable and the pair $(A,Q^{1/2})$ is detectable. In this case, for a given state $y$, the feedback control is determined by the following relation
\begin{equation}
   u(y) = -R^{-1} B^{\top} Py,
\label{control_SDRE}
\end{equation}
where $P\in\mathbb{R}^{N_h\times N_h}$ is the unique positive semidefinite solution of the Algebraic Riccati Equation (ARE)
\begin{align*}
A^\top P + P A
  -PBR^{-1}B^\top P + Q = 0.
\end{align*}
For the sake of notation, let $F=BR^{-1}B^\top$.
In essence, the SDRE technique extends this approach by incorporating a state-dependence in the Riccati solution, leading to the following feedback law:
\begin{equation}
u(y) = -R^{-1} B^{\top} P(y)y,
\label{control_sdre_feed}
\end{equation}
where $P(y)$ now serves as the solution to a 
{SDRE defined as}
\begin{align}
A^\top(y) P(y) + P(y) A(y)-P(y)F P(y)& + Q = 0,
\label{sdre}
\end{align}
with
$A(y)$ fixed at the state $y$. 
In general, this iterative process can be applied along the {parametric} trajectory, sequentially solving \eqref{sdre} as the state $y_{\mu,t}$ evolves. The controlled dynamics now reads
\begin{equation}
    \label{semilinear_feed}
 \begin{cases}
 \dot{y}(t; \mu) & = A_{cl}(y(t;\mu))y(t;\mu), \\
    \second{y(0;\mu)} & =\second{ y_{\mu,0}},    
 \end{cases}
\end{equation}
where $A_{cl}(y(t; \mu)) = A(y(t;\mu))  -F P(y(t;\mu))$ is called closed-loop matrix.
We assume that the stabilizability and detectability conditions are satisfied for all states $y$ within a domain $\Omega$.
\begin{assump}
\label{assum:stab_detec}
    The pair $(A(x),B)$ and $(A(x),Q^{1/2})$ are \first{stabilizable and detectable, respectively, for every $x \in \Omega$}.
\end{assump}
Under these conditions, the SDRE \eqref{sdre} is guaranteed to have a unique positive semidefinite solution for all $y \in \Omega$ (see, e.g., \cite{lancaster1995algebraic}).
Under the additional assumption that $f \in C^1(\mathbb{R}^{N_h})$, it can be shown that the controlled dynamics arising from the system in \eqref{semilinear_feed} converges locally to the target state (for precise statements and further elaboration, see \cite{Banks_Lewis_Tran_2007}).

\first{It can be shown that the closed-loop system governed by the feedback law \eqref{control_sdre_feed} exhibits local asymptotic stability, as formalized in the following theorem from \cite{Banks_Lewis_Tran_2007}.

\begin{thm}
    Assume that the system \eqref{eq:dynamic} is such that \( f(x) \) and \( \frac{\partial f(x)}{\partial x_j} \) (for \( j = 1, \dots, n \)) are continuous in \( x \) for all \( \|x\| \leq \hat{r} \), where \( \hat{r} > 0 \).
Assume further that \( A(x) \) and \( B(x) \) are continuous in some nonempty neighborhood of the origin \( \mathcal{O} \subseteq B_{\hat{r}}(0) \) and that Assumption \ref{assum:stab_detec} holds.. Then, the system with the control given by \eqref{control_sdre_feed} is locally asymptotically stable.
\label{thm:stable}
\end{thm}

Global stabilization can be achieved by introducing additional assumptions. For example, this holds if the closed-loop matrix is symmetric or diagonally dominant with all diagonal elements negative for all $x\in \mathbb{R}^{N_h}$. We refer the reader to \cite{cloutier1996nonlinear,chen2015global} for further discussion on the topic and additional cases.
}

We now emphasize the dependence of the Riccati solution on the parameter $\mu$ and time $t$, denoting it as 
$P_{\mu,t} = P(y(t;\mu))$. 
This distinction is crucial for the findings in this contribution, particularly for the tailored approaches to solve parametric nonlinear feedback control problems by employing SDREs. 
In what follows, we describe the numerical solver we used to deal with the solution of \eqref{sdre}.

\subsection{Numerical approximation of the SDRE}
\label{sec:NK}
Motivated by the parametric and time-dependent nature of our setting, which results in a repeated high computational effort, we propose to solve the SDRE in \eqref{sdre} using an iterative method instead of a direct one. 
Specifically, we pay peculiar attention to the choice of the \emph{initial guess} of the nonlinear solver. Indeed, the Riccati solution approximated at a given time instance and parameter value can serve as an effective initial guess for approximating the Riccati solution at a different time or parameter value. This methodology will be further clarified in Section \ref{sec:sensitivity}, where we establish conditions under which such initialization enhances the effectiveness of the iterative scheme. 
First of all, let us define the \emph{residual} of the SDRE as 
\begin{align}
\mathcal R_{(\mu, t)}(P_{\mu, t}) = A^\top_{\mu,t} P_{\mu, t} + P_{\mu, t} A_{\mu,t}-P_{\mu, t}F P_{\mu, t} + Q,
\label{eq:ResSDRE}
\end{align}
with $A_{\mu,t}=A(y_{\mu,t})$.
In other words, given a pair $(\mu,t)  \in \mathcal P \times \mathbb{R}^+$, solving the SDRE translates into finding the solution to $\mathcal R_{(\mu, t)}(P_{\mu, t}) = 0$. Due to the nonlinear nature of \eqref{eq:ResSDRE}, a Newton method can be applied to reach an approximation of the solution. In particular, given an initial guess $P^0_{\mu,t}$, at the $(k+1)$-th iteration of the method we solve 
\begin{equation}
    \label{eq:Newton}
    \mathcal R_{(\mu, t)}'[P_{\mu, t}^k](\delta P_{\mu, t}) = - \mathcal R_{(\mu, t)}(P_{\mu, t}^k) \quad \text{with the update}\quad P_{\mu, t}^{k+1} = P_{\mu, t}^k +  \delta P_{\mu, t},
\end{equation}
where $\mathcal R_{(\mu, t)}'[P_{\mu, t}]$ denotes the Fr\'echet derivative of the residual.
Specifically, at each iteration, given the previous iterate $P^k_{\mu, t}$, we solve the Lyapunov equation
\begin{equation}
    \label{eq:FrechetSDRE}
    A_{cl}^{k}(y_{\mu,t})^\top P^{k+1}_{\mu, t} + P^{k+1}_{\mu, t}A_{cl}^{k}(y_{\mu,t}) = -P^{k}_{\mu, t} F P^{k}_{\mu, t} -Q,
\end{equation}
with $ A_{cl}^{k}(y_{\mu,t}) = A_{\mu,t}-F P^{k}_{\mu, t}$,
generating a sequence $\{P^k_{\mu, t}\}_{k\geq 0}$.
It is possible to prove that the sequence $\{P^k_{\mu, t}\}_{k\geq 0}$ converges to the unique solution of the SDRE (see, e.g., \cite{kleinman1968iterative}), \first{under the following assumptions.}
\begin{thm}
\label{thm:NKconvergence}
Let $P^0_{\mu, t} \in \mathbb{R}^{N_h \times N_h}$ be symmetric and positive semideﬁnite and such that $ A_{\mu,t}-F P^{0}_{\mu, t}$ is stable, and let Assumption \ref{assum:stab_detec} hold. Then the Newton iterates $\{P^k_{\mu, t}\}_k$ \first{converge} to the unique solution of the SDRE. 
\end{thm}

 This strategy, denoted as Newton-Kleinman (\texttt{NK}) \cite{kleinman1968iterative}, is summarized in Algorithm \ref{alg:newton}. 
 We remark that the Lyapunov equation is solved via the Matlab function \texttt{lyap}, based on the Bartels–Stewart's strategy \cite{bartels1972algorithm}.

\begin{algorithm}
\caption{The \texttt{NK} algorithm}
\label{alg:newton}
\begin{algorithmic}[1]
\Statex{\textbf{Input:} $A(\cdot)$, $F$, $Q$, tolerance $\mathsf {tol_{\text{\texttt{NK}}}}$, initial guess $P^0_{\mu,t}$} 
\Statex{\textbf{Output:} solution of the SDRE, i.e., $P_{\mu, t}$}
\Statex
\State{\first{$j=0$}}
\While{$\norm{\mathcal R_{\mu,t}(\first{P^{j}_{\mu, t}})} \geq \mathsf {tol_{\text{\texttt{NK}}}}$}
\State{Given $\first{P^{j}_{\mu, t}}$, solve \eqref{eq:FrechetSDRE} to get $\first{P^{j+1}_{\mu, t}}$}
\State{\first{$j = j + 1$}}
\EndWhile
\end{algorithmic}
\end{algorithm} 
    As in any Newton-based method, the initial guess plays a fundamental role in the \texttt{NK} strategy. Indeed, a proper choice for $P^0_{\mu,t}$ in the \texttt{NK} algorithm can significantly accelerate the convergence towards the unique positive semidefinite solution of the SDRE. In this work, we propose a \emph{parameter-based selection of the guess}. \both{For the sake of presentation, we restrict to one-dimensional parameter spaces. This choice is not restrictive: the concept is suited and generalizable to the multidimensional parametric space as done, for example, in the numerical results of Section \ref{sec:results}.} Given two different parameters $\mu_i$ and $\mu_{i+1}$ with $\mu_i < \mu_{i+1}$, we solve the SDRE for $P_{\mu_{i+1},t}$ by leveraging the information from $P_{\mu_{i},t}$ for a given $t$.
    Following the ideas behind \emph{continuation} methods \cite{KellerLecturesNumericalMethods1988}, we propose using $P_{\mu_{i},t}$ as the initial guess for solving the SDRE for $P_{\mu_{i+1},t}$, thereby ensuring faster convergence of the \texttt{NK} towards the next parametric solution of the SDRE. Moreover, in Section \ref{sec:sensitivity}, we demonstrate that, {under suitable assumptions}, employing $P_{\mu_{i},t}$ as the initial guess ensures the stability of the initial closed-loop matrix $A_{cl}^{0}(y_{\mu_{i+1},t})$, thereby guaranteeing the convergence of the scheme. 
    We refer to this approach as the \emph{cascade-Newton-Kleinman} (\texttt{C-NK}) method,  which is summarized in Algorithm \ref{alg:cnk}.
    This concept can be extended to consecutive time instances for a fixed parameter $\mu$. Specifically, given $t_{i} < t_{i+1}$, we can use $P_{\mu,t_i}$ as initial guess when solving for $P_{\mu,t_{i+1}}$.
    In Section \ref{sec:sensitivity}, theoretical justifications for selecting the distance between $\mu_i$ and $\mu_{i+1}$ (or $t_i$ and $t_{i+1}$) are provided. 
    This section explores technical aspects that, while not essential for all readers, provide valuable insights into the underlying methodology. Readers primarily interested in the main DLRA-based results and their implications may proceed directly to Section \ref{sec:DOrb} without any loss of continuity.

\begin{algorithm}[H]
\caption{The \texttt{C-NK} algorithm}
\label{alg:cnk}
\begin{algorithmic}[1]
\Statex{\textbf{Input:} $A(\cdot)$, $F$, $Q$, tolerance $\mathsf {tol_{\text{\texttt{NK}}}}$, initial guess $P^0_{\mu_1,t}$, \underline{ordered} set $\mathcal P_h = \{\mu_k\}_{k=1}^p$}
\Statex{\textbf{Output:}  $P_{\mu,t}$ for all $\mu \in \mathcal P_h$}
\Statex{}
 \For{$k = 1, \ldots, p$}
    \State{$P_{\mu_{k},t}$ = \texttt{NK}$(A(\cdot), F, Q, \mathsf {tol_{\text{\texttt{NK}}}}, P^0_{\mu_k,t})$}
    	\State{$P^0_{\mu_{k+1},t} = P_{\mu_k,{t}}$}    
\EndFor
\end{algorithmic}
\end{algorithm}

    

    


\section{Sensitivity estimates}
\label{sec:sensitivity}

In this section, we derive sensitivity estimates within the SDRE framework to provide a theoretical foundation for the cascade approach. For enhanced readability, all auxiliary results required for these estimates \both{and} their numerical validation are presented in Appendix \ref{sec:appendix}. Section \ref{sec:gen_estimates} outlines the definitions, assumptions, and established results that form the basis for deriving the sensitivity estimates of interest. In Section \ref{sec:param}, we propose the sensitivity estimates related to the employment of the cascade information in the parametric space, while Section \ref{sec:time} deals with the cascade estimates for time evolution. 
\subsection{Preliminary notions}
\label{sec:gen_estimates}
\second{In the following, we provide some preliminary notions that will be necessary to the proposed sensitivity estimates.}

\begin{definition}
    
A vector field \( f(t, y) \) satisfies a one-sided Lipschitz condition in a domain $\Omega \subset \mathbb{R}^d$ if there exists a constant \( L \in \mathbb{R} \) such that
\begin{equation}
 (f(t, y_1) - f(t, y_2)) \cdot (y_1 - y_2) \leq L \| y_1 - y_2 \|^2, \quad \forall t \in \mathbb{R}^+, \; \forall y_1, y_2 \in \Omega .
 \label{eq:1side_lip}
\end{equation}
\end{definition}

Notably, the one-sided Lipschitz constant $L$ differs from the standard Lipschitz constant since it can be negative. This distinction is fundamental, as the sign of $L$ plays a key role in determining whether the system is dissipative, which is expected in a controlled problem.
\first{In the linear case, the notion of one-sided Lipschitz constant is equivalent to the logarithmic norm of $A$}, introduced in the next definition.

\begin{definition}
    Given a matrix $A\in \mathbb{C}^{d \times d}$, the logarithmic norm of the matrix $A$ is \both{defined as}
\begin{equation}
\label{log_norm}
\both{\textup{logm}}(A)= \sup_{x \in \mathbb{C}^d \setminus \{0\}} \frac{\both{\textup{Re}} \left( ( Ax) \cdot x \right)} {\Vert x \Vert^2},
\end{equation}
\end{definition}
where $\both{\textup{Re}}(z)$ indicates the real part of $z \in \mathbb{C}$.

The logarithmic norm is a crucial tool in the stability analysis of both continuous and discrete linear dynamical systems. Specifically, since it can be demonstrated that $\Vert e^{tA} \Vert \le e^{t \,\both{\textup{logm}}(A)}$ for all $t \ge 0$ (see, for example, \cite{soderlind2006logarithmic}), we can assert that the dynamical system is stable if $\both{\textup{logm}}(A) \le 0$.
\second{Throughout this section, we consider the following assumptions.
\begin{assump}
\color{white}{text} \color{black}\\
\label{assum:review}
   \begin{enumerate}[(a)]
  
    \item The matrix function $A_{cl}(y) = A(y)-F P(y)$ is Lipschitz continuous in $\Omega \subset \mathbb{R}^d$  with Lipschitz constant $L_{cl}>0$, i.e.,
    \begin{equation}
        \Vert A_{cl}(x)-A_{cl}(y) \Vert \le L_{cl} \Vert x -y \Vert.
    \end{equation}
    \label{hp:Acl_lip}
    \item The pair $(A(x),B)$ is stabilizible for every $x \in \Omega \subset \mathbb{R}^d$.
    \label{hp:stab}
    \item The pair $(A(x),Q^{1/2})$ is detectable for every $x \in \Omega \subset \mathbb{R}^d$.
    \label{hp:detec}
    \item The function $y_0(\mu): \mathcal{P} \subset \mathbb{R}^p \rightarrow \mathbb{R}^d$ is Lipschitz continuous with Lipschitz constant $L_{y_0}>0$.
    \label{hp:y0_lip}
    \item The matrix function $A(y)$ is Lipschitz continuous in $\Omega \subset \mathbb{R}^d$ with Lipschitz constant $L_A>0$,  \label{hp:A_lip}
    \item\both{Given $(\nabla A_{cl}(y)) y =
\begin{pmatrix}
       \left [\frac{\partial}{\partial y_1}A_{cl}(y) \right ] y
    \ldots
    \left [\frac{\partial}{\partial y_d}A_{cl}(y) \right ] y
  \end{pmatrix} \in \mathbb{R}^{d \times d}
$ for $y \in \mathbb{R}^{d}$}, we assume
\begin{equation*}
\both{\textup{logm}}(A_{cl}(y))+\both{\textup{logm}}((\nabla A_{cl}(y)) y) <0, \quad \forall y \in \Omega \subset \mathbb{R}^d.
\end{equation*}

\label{hp:nabla}
\end{enumerate}
\end{assump}
}

\noindent \second{Assumption 2 \ref{hp:stab} and \ref{hp:detec} ensure the existence and uniqueness of the solution to the SDRE. Consequently, they are regarded as standing assumptions throughout this section.}

 For the estimates, we consider the solutions of two SDREs, defined by the pairs $(\mu,t)$ and $(\tilde{\mu},\tilde{t})$ as
\begin{equation}
A^\top_{\mu,t}  P_{\mu,t} + P_{\mu,t} A_{\mu,t}-P_{\mu,t}F P_{\mu,t}  + Q = 0,
\label{riccati_mu1}
\end{equation}
and 
\begin{equation}
A^\top_{\tilde{\mu},\tilde{t}} P_{\tilde{\mu},\tilde{t}} + P_{\tilde{\mu},\tilde{t}} A_{\tilde{\mu},\tilde{t}}-P_{\tilde{\mu},\tilde{t}}F P_{\tilde{\mu},\tilde{t}} + Q = 0,
\label{riccati_mu2}
\end{equation}
where $A_{\mu,t} = A(y(t;\mu))$ and $P_{\mu,t} = P(y(t;\mu))$, using the same notation introduced in Section \ref{sec:problem}. 
For both estimates, in the parametric space and time domains, we utilize a foundational result derived from Lemma 3.1 in \cite{gahinet1990computable} and Theorem 2.4 in \cite{kenney1990sensitivity}. Indeed, denoting the perturbation $A_{\mu,t}-A_{\tilde{\mu},\tilde{t}}$ as $\Delta A_{\mu,\tilde{\mu},t, \tilde{t}}$, the following theorem holds.

\begin{thm}
\label{thm:deltaA}
Suppose \eqref{riccati_mu1} has a unique positive semidefinite solution $P_{\mu,t}$. Let us assume that the perturbation
$\Delta A_{\mu,\tilde{\mu},t,\tilde{t}}$  is such that
\second{
\begin{equation}
    \Vert \Delta A_{\mu,\tilde{\mu},t,\tilde{t}} \Vert < \frac{1}{2 \Vert H_{\mu,t} \Vert},
    \label{condi_pert}
\end{equation}
}
where $H_{\mu,t}$ is the solution
of the following Lyapunov equation

\begin{equation}
   \second{\Omega_{P_{\mu,t}}(H_{\mu,t})=} (A_{\mu,t}- FP_{\mu,t})^\top H_{\mu,t} + H_{\mu,t}(A_{\mu,t}- FP_{\mu,t}) = - I,
    \label{eq:lyap_pert}
\end{equation}
\second{and $\Vert \cdot \Vert$ denotes the operator 2-norm.}
Then \eqref{riccati_mu2} has a unique symmetric \both{positive semidefinite} solution $P_{\tilde{\mu},\tilde{t}}$ and $A_{\tilde{\mu},\tilde{t}}- F P_{\mu,t}$ is stable.
\end{thm}

\second{
 \begin{proof}
     By Lemma~3.1 in~\cite{gahinet1990computable}, in the absence of perturbations in the matrix \( F \) (i.e., \( \Delta F = 0 \)), under the condition that
     \begin{equation}
    \Vert \Delta A_{\mu,\tilde{\mu},t,\tilde{t}}  \Vert < \frac{1}{2  \Vert \Omega_{P_{\mu,t}}^{-1}  \Vert},
    \label{condi_pert2}
\end{equation}
where the operator norm of \( \Omega_{P}^{-1} \) is defined by
\[
    \Vert \Omega_P^{-1}  \Vert
    = \sup_{\substack{X \in \mathbb{R}^{d \times d} \\ X \neq 0}} 
       \frac{ \Vert \Omega_P^{-1}(X) \Vert}{ \Vert X  \Vert},
\]
the perturbed Riccati equation~\eqref{riccati_mu2} admits a unique symmetric positive semidefinite solution and the matrix \( A_{\tilde{\mu},\tilde{t}} - F P_{\mu,t} \) is stable.
Moreover, by Theorem~2.4 in~\cite{kenney1990sensitivity}, the operator norm \(  \Vert \Omega_{P_{\mu,t}}^{-1}  \Vert \) coincides with the spectral norm of the unique solution \( H_{\mu,t} \) to the Lyapunov equation~\eqref{eq:lyap_pert}. Hence, condition~\eqref{condi_pert2} is equivalent to condition~\eqref{condi_pert}, which completes the proof.

 \end{proof}
     }

In the next section, we aim to demonstrate that, under specific conditions, the Riccati equation approximation for a given parameter $\mu$ also serves as a stabilizing solution for a different parameter $\tilde{\mu}$, for a fixed time.

\subsection{Cascade sensitivity estimates in the parameter space}
\label{sec:param}
In the \texttt{C-NK} approach, we first approximate the SDRE for the first parameter of a finite parametric set $\mathcal P_h \subset \mathcal P$ and then use it as an initial guess for the next SDRE approximation via the \texttt{NK} method. Suppose that $\mathcal{P}_h$ is obtained by discretizing $\mathcal P$ with a uniform stepsize $\Delta \mu$.
We now apply the concepts introduced in Section \ref{sec:gen_estimates} to the \texttt{C-NK} framework, deriving a criterion for selecting $\Delta \mu$ that ensures the stabilizability of the closed-loop matrix at the initial time. Furthermore, under a suitable assumption on the one-sided Lipschitz constant, this result can be extended to the entire time interval, thereby guaranteeing the convergence of the \texttt{C-NK} algorithm over time. For simplicity, we restrict our focus to one-dimensional parameter spaces. However, the proofs readily generalize to multidimensional cases.

Now, we can prove the following Proposition, which provides the justification \first{for} the cascade information in the \texttt{NK} procedure. The auxiliary results needed for the proof are provided in Appendix \ref{sec:appendix}.
\begin{proposition}
\label{prop:dmu}
\second{Given Assumption  \ref{assum:review} \ref{hp:Acl_lip}-\ref{hp:A_lip}}, at the initial time, given $\mu \in \mathcal P$ and $\tilde{\mu} = \mu + \Delta \mu$, if the parameter stepsize is chosen such that
\begin{equation}
\Delta \mu \second{<} \frac{ 1}{2 L_A L_{y_0}\Vert\second{ H_{\mu,0} }\Vert},
\label{cond_init}
\end{equation}
then $A_{\tilde{\mu},0}-F P_{\mu,0}$ is stable.
Moreover, if 

\begin{equation}
L  \second{<} - \frac{1}{t}\both{\log}(2 L_A  L_{y_0} \Delta \mu \Vert {H_{\mu,t}} \Vert), \quad t \in (0,T],
\label{cond_L}
\end{equation}
then $A_{\tilde{\mu},t}-B P_{\mu,t}$ is stable for all $t\in[0,T]$.

Finally, provided that \eqref{cond_init} and \eqref{cond_L} hold, the Newton-Kleinman sequence $\{P^j_{\tilde{\mu},t}\}_j$ arising from the resolution of \eqref{eq:FrechetSDRE} with initial guess $P_{\mu,t}$ converges to the unique solution of the SDRE \eqref{riccati_mu2} for all $t \in [0,T]$.
    \end{proposition}

\begin{proof}
By Assumption \ref{assum:review} \ref{hp:A_lip} and \second{Corollary \ref{cor:decay_lip_1} in Appendix \ref{sec:ar}} we have
$$
    \Vert \Delta A_{\mu,\tilde{\mu},t,t} \Vert \le L_A \Vert y(t;\mu)-y(t;\tilde{\mu}) \Vert \le L_A  L_{y_0} \Vert \mu - \tilde{\mu} \Vert e^{Lt} = L_A  L_{y_0} \Delta \mu \, e^{Lt} .
    $$

Condition \eqref{condi_pert} is satisfied if 
\begin{equation}
L_A  L_{y_0} \Delta \mu \, e^{ Lt} \second{<} \frac{1}{2 \Vert H_{\mu,t} \Vert}.
\label{cond_proof}
\end{equation}
Condition \eqref{cond_init} is obtained by imposing $t=0$ in \eqref{cond_proof},
while condition \eqref{cond_L} is obtained manipulating \eqref{cond_proof} for $t>0$.
The convergence of the Newton-Kleinman sequence follows from Theorem \ref{thm:NKconvergence}, under assumptions \ref{hp:stab} and \ref{hp:detec}, with the additional observation that the initial guess $P_{\mu,t}$ is such that  $A_{\tilde{\mu},t}- F P_{\mu,t}$ is stable.
\end{proof}
\begin{remark}
\label{rem:rem_param}
\first{From a numerical perspective, the dynamical system is integrated in time using an appropriate quadrature scheme.}
     Similar results can be obtained substituting the solution $y(t;\mu)$ with an approximation $\tilde{y}(t;\mu)$ of order $p$ \first{of accuracy in time} satisfying $\Vert y(t;\mu)- \tilde{y}(t;\mu)\Vert \le C (\Delta t)^p$. Indeed, since 
     \begin{equation}
     \label{eq:rem1_1}
   \Vert \tilde{y}(t;\mu) - \tilde{y}(\tilde{t}; \tilde{\mu}) \Vert \le  \Vert y(t;\mu) - y(\tilde{t}; \tilde{\mu}) \Vert  + 2C (\Delta t)^p,
     \end{equation}
     the matrix perturbation can be bounded as
     \begin{equation}
     \label{eq:rem1_2}
     \Vert \Delta A_{\mu,\tilde{\mu},t,t} \Vert \le L_A \Vert \tilde{y}(t,\mu)-\tilde{y}(t,\tilde{\mu}) \Vert \le L_A ( L_{y_0} \Vert \mu - \tilde{\mu} \Vert e^{Lt} + 2C (\Delta t)^p),
    \end{equation}
    and thus, condition \eqref{cond_L} becomes
\begin{equation}
L  \second{<} \frac{1}{t}[\both{\log}(1-\second{4}C(\Delta t)^p L_A \Vert H_{\mu,t} \Vert) - \both{\log}(2 L_A  L_{y_0} \Delta \mu \Vert H_{\mu,t} \Vert)], \quad t \in (0,T].
\label{cond_L_discrete}
\end{equation}
     
 \end{remark}
\subsection{Cascade sensitivity estimates in time}
\label{sec:time}
The sensitivity estimates over the parametric space can also be extended to the time variable. We want to prove that the solution $P_{\mu,t}$ to \eqref{riccati_mu1} at time instance $t$ is stabilizing for the Riccati equation \eqref{riccati_mu2} with $\tilde{t} = t + \Delta t$, for a proper timestep $\Delta t$ and a fixed parameter $\mu$ (i.e., $\mu = \tilde{\mu}$).
We aim to apply Theorem \ref{thm:deltaA} to derive a criterion for appropriately selecting the timestep $\Delta t$ in order to guarantee a stabilizing guess for the \texttt{NK} algorithm in time. This result is formalized in Proposition \ref{prop:Delta_t}.
The auxiliary results needed to prove this statement are detailed in Appendix \ref{sec:appendix}.

\begin{proposition}
\label{prop:Delta_t}
\second{Given Assumption \ref{assum:review} \ref{hp:stab}, \ref{hp:detec}, \ref{hp:A_lip} and \ref{hp:nabla}},
if the timestep $\Delta t$ for a given parameter $\mu$ is selected to satisfy the condition
\begin{equation}
\Delta t \second{<} \frac{1}{2 L_A \Vert A_{cl}(y_{\mu,0})y_{\mu,0} \Vert \Vert H_{\mu,t} \Vert },
\label{cond_init_dt}
\end{equation}
then  $A_{\mu,t+\Delta t}-F P_{\mu,t}$ is stable.


Finally, provided that \eqref{cond_init_dt} holds, the Newton-Kleinman sequence $\{P^j_{\mu,t+\Delta t}\}_j$ arising from the resolution of \eqref{eq:FrechetSDRE} with initial guess $P_{\mu,t}$ converges to the unique solution of the SDRE \eqref{riccati_mu2}.
    \end{proposition}

\begin{proof}
   We define $\Delta y(t) = y_{\mu, t+\Delta t}-y_{\mu, t}$.
   We know that it holds
   $$
   \Vert \Delta y(t) \Vert \le e^{L t} \Vert \Delta y(0) \Vert =  e^{L t} \Vert y_{\mu,\Delta t} - y_{\mu,0} \Vert \le e^{L t} \int_0^{\Delta t} \Vert A_{cl} (y_{\mu,s})y_{\mu,s} \Vert ds \le \Delta t e^{L t} M_{\mu, \Delta t},
   $$
where $ M_{\mu, \Delta t} =  \max_{\overline{t} \in [0, \Delta t]} \Vert A_{cl}(y(\overline{t},\mu))y(\overline{t},\mu) \Vert $. \second{By Proposition \ref{prop:norm_Ay} in Appendix \ref{sec:ar}} we have $ M_{\mu, \Delta t} = \Vert A_{cl}(y_{\mu,0})y_{\mu,0} \Vert$, while by Assumption \ref{assum:review} \ref{hp:A_lip} we have
$$
    \Vert \Delta A_{\mu,\mu,t,\tilde{t}} \Vert \le L_A \Vert \Delta y(t) \Vert \le  L_A \Delta t e^{L t} M_{\mu, \Delta t} ,
    $$
    then condition \eqref{condi_pert} is satisfied if condition \eqref{cond_init_dt} holds.
    The proof of the convergence of the Newton-Kleinman sequence relies on similar arguments as those presented in Proposition \ref{thm:deltaA} and is therefore omitted.
\end{proof}

 \begin{remark}
     Similar results can be obtained by substituting the solution $y(t;\mu)$ with an approximation $\tilde{y}(t;\mu)$ of order $p$ \first{of accuracy in time}, using the same arguments of Remark \ref{rem:rem_param}. Indeed, assuming \eqref{eq:rem1_1}, when selecting the timestep $\Delta t$, we obtain
      
$$
    \Vert \Delta A_{\mu,\mu,t,\tilde{t}} \Vert \le  L_A \Delta t ( e^{L t} M_{\mu, \Delta t} +  2C (\Delta t)^{p-1}),
$$
then, under the assumptions $\Delta t <1$ and $p \ge 1$, condition \eqref{cond_init_dt} becomes 
$$
\Delta t \second{<} \frac{1}{2 L_A (M_{\mu, \Delta t} + 2C )\Vert H_{\mu,t} \Vert }.
$$
     
 \end{remark}

\section{Dynamical Low Rank Approximation for feedback control}
\label{sec:DOrb}
This section focuses on the application of the DLRA in the feedback control setting. We start from the basics of this approach in Section \ref{sec:DLRA_basic}, and, then, we specialize the results to semilinear parametric infinite horizon optimal control problems in Section \ref{sec:DLRA_ex}. 
\subsection{System decomposition and tangent spaces}
\label{sec:DLRA_basic}
Let us suppose we want to solve the parametric system \eqref{eq:dynamic}, for a fixed set of parameters $\{\mu_j\}_{j=1}^p \subset \mathcal P$, for $j$ and $p$ in $\mathbb N$.
For \emph{a sampled} initial condition matrix $\mathcal Y_0(\mu) = [y_0^{\mu_1}, \dots, y_0^{\mu_p}] \in \mathbb R^{N_h \times p}$, we aim at finding $\mathcal Y = \mathcal Y(\mu) \in \mathcal C^1((0, \infty), \mathbb R^{N_h \times p})$ such that

\begin{equation}
\label{eq:DLRA_dynamic}
\begin{cases}
    \dot{\mathcal Y}(t) = \mathcal F(\mathcal Y(t))  + B \mathcal{U} (\mathcal Y(t)) & \text{for } t \in (0, \infty), \\
    \mathcal Y(0) = \mathcal Y_0(\mu),
\end{cases} 
\end{equation}
where $\mathcal U: \mathbb R^{N_h \times p}  \rightarrow \mathbb R^{m \times p}$ is the \emph{control matrix},  $B \in \mathbb R^{N_h \times m}$ and $\mathcal F (\mathcal Y(t)) =\left[ f(y(t;\mu_1)), \ldots , f(y(t;\mu_p)) \right] \in  \mathbb R^{N_h \times p}$. 
In other terms, we recast the parametric problem \eqref{eq:dynamic}, considering the dynamics for a set of $p$ parameters \emph{simultaneously}.
The main objective is to represent the high-dimensional solution $\mathcal Y(t)$ in a reduced space framework of dimension $r \ll N_h$ as follows: 
\begin{equation}
\label{eq:UZt}
    \mathcal Y(t) \approx \mathsf Y = \mathsf  U\mathsf Z^{\top} = \sum_{i=1}^r U_i(t)Z_i(t;\mu),
\end{equation}
where \first{$\mathsf U = [U_1(t)|\dots| U_r(t)] \in \mathcal C^1((0, \infty),\mathbb R^{N_h \times r})$} and \first{$\mathsf Z \in \mathcal C^1((0, \infty),\mathbb R^{p\times r})$} with $\mathsf Z_{ji}(t) = Z_i(t; \mu_j)$. \first{For the sake of notation, we are dropping the time and the parameter dependencies from the matrices $\mathsf Y$, $\mathsf U$, and $\mathsf Z$.} Furthermore, $\mathsf U$ has orthonormal columns, i.e., $\mathsf U^\top \mathsf U = \mathsf I_r$, where $\mathsf I_r$ is the identity matrix of dimension $r$.
In other words, we represent the reduced solution in the manifold
\begin{equation}
\label{eq:Mr}
    \mathcal M_r = \{ \mathsf Y \in \mathbb R^{N_h\times p} \; : \; \mathsf Y = 
    \mathsf  U \mathsf  Z^\top, \; \mathsf U^\top \mathsf U = \mathsf I_r, \; \mathsf  U \in \mathbb R^{N_h\times r}, \;  \mathsf  Z \in \mathbb R^{p\times r}\}.
\end{equation}
\first{To avoid ambiguity, we emphasize that $\mathsf U$ refers to a basis function matrix and not to the control action, which is denoted with $\mathcal U$.}
The \first{approximation} \eqref{eq:UZt} is unique up to \first{rotations and reflections} via an orthonormal matrix (see, e.g., \cite{KocLubich}).
Thus, as an alternative to the non-unique decomposition \eqref{eq:UZt}, the DLRA
looks for a unique representation of the \emph{tangent space} of \eqref{eq:Mr} at $\mathsf Y = \mathsf U \mathsf Z^{\top}$ \cite{KocLubich,sapsis2009dynamically}, defined as
\begin{equation}
    \label{eq:TMn}
    T_{\mathsf Y} \mathcal M_r 
    = \{
    \dot{\mathsf Y} \in \mathbb R^{N_h \times p} \; : \; \dot{\mathsf Y} = \dot{\mathsf U} \mathsf Z^{\top} + \mathsf U \dot{\mathsf Z}^{\top},  \; \dot{\mathsf U} \in \mathbb R^{N_h \times r}, \; \dot{\mathsf Z} \in \mathbb R^{p \times r}, \; \mathsf U^{\top} \dot{\mathsf U} \in \text{so}(r)
    \}, 
    \end{equation}
where $\text{so}(r)$ is the space of the skew-symmetric real matrices in $\mathbb{R}^{r\times r}$. 

However, the definition of this tangent space is still not unique. 
To guarantee a unique representation of the tangent space \eqref{eq:TMn}, we further assume an orthogonality constraint on $\mathsf U$, i.e., $\mathsf U^{\top} \dot{\mathsf U} = 0$.
Thanks to this further assumption, with simple algebraic manipulation, we can provide the following differential equations for $\mathsf U$ and $\mathsf Z$ \cite{lubich06,Pagliantini2021409}:
\begin{equation}
\label{eq:UZ_dynamics}
\begin{cases}
    \dot{\mathsf U} = (\mathsf I_{N_h}- \mathsf U \mathsf U^{\top}) \dot{\mathsf Y} \mathsf Z \mathsf C^{-1}, & \\
    \dot{\mathsf Z} = \dot{\mathsf Y}^{\top} \mathsf U, &
\end{cases}   
\end{equation}
where $\mathsf C = \mathsf Z^{\top} \mathsf Z$. \\
\first{The use of the tangent space to represent $\dot{\mathsf Y}$ is motivated by some intriguing properties such as \cite{KocLubich,lubich06,Pagliantini2021409} (i) the computational feasibility: the algorithm avoids costly full matrix operations, requiring only multiplications with thin matrices and a low-rank decomposition at the initial time, and (ii) natural applicability to dynamical systems: when $\mathcal Y(t)$ is the solution to a differential equation $\dot{\mathcal Y} = F(\mathcal Y)$ one approximates the dynamics instead of $\dot{\mathcal Y}(t)$, enabling model reduction directly within the evolution law and increasing the accuracy of the reduced solution, as we will illustrate in Section \ref{sec:results}.}
\subsection{Differential equations for the reduced representation}
\label{sec:DLRA_ex}
In this section, we formalize the reduction of the infinite horizon feedback control problem, as introduced in Section \ref{sec:problem}, using the DLRA methodology.\\
The aim is to provide an initial decomposition for $\mathsf Y_0 = \mathsf U_0 \mathsf Z^\top_0$ and to evolve it in time. First, we perform an SVD on the initial datum $\mathcal Y_0(\mu)$. The matrix $\mathsf U_0$ is the matrix formed by the first $r$ left singular vectors of $\mathcal Y_0(\mu)$, while the coefficient matrix is given by $\mathsf Z_0 = \mathcal Y_0^{\top} \mathsf U_0$.\\
In our setting, we exploit the definition of the system \eqref{eq:DLRA_dynamic} approximated by
\begin{equation}
\label{eq:approx_dyn}
\begin{cases}
    \dot{\mathsf Y} = \mathcal F (\mathsf Y)  + B \mathcal U(\mathsf Y),\\
    \mathsf Y(0) = \mathsf Y_0.
\end{cases}
\end{equation}
Substituting the reduced dynamics \eqref{eq:approx_dyn} into equation \eqref{eq:UZ_dynamics}, we obtain
\begin{equation}
\label{eq:UZ_dynamics_feed}
\begin{cases}
    \dot{\mathsf U} = (\mathsf I_{N_h} - \mathsf U \mathsf U^T) 
    ( 
    \mathcal F(\mathsf U \mathsf Z^{\top})  + B \mathcal U (\mathsf U \mathsf Z^{\top})) \mathsf Z \mathsf C^{-1},  & \\
 \dot{\mathsf Z} =  ( \mathcal F(\mathsf U \mathsf Z^{\top})+ B\mathcal U(\mathsf U\mathsf Z^{\top}))^\top \mathsf U. &
\end{cases}   
\end{equation}
We now want (i) to derive a feedback law for the control matrix $\mathcal U (\cdot)$ and (ii) to establish a specific low rank and efficient representation of \eqref{eq:UZ_dynamics_feed} when $f (y)= A y + T \left( y \otimes y \right) $, where $A \in \mathbb{R}^{N_h \times N_h}$ and $T \in \mathbb{R}^{N_h \times N^2_h}$.

This is the case of the numerical results presented in Section \ref{sec:results}. Nevertheless, this approach can be readily extended to polynomial nonlinearities; however, we note that quadratic nonlinearities are prevalent in many applications. 

Let us start with (i) the definition of the feedback law.
By transposing the equation for $ \mathsf Z$ and denoting $\widetilde{\mathsf Z}= \mathsf Z^\top$, we obtain for each parameter $\mu_j$
$$
\dot{\widetilde{\mathsf Z}}_j = \mathsf U^\top A \mathsf U \widetilde{\mathsf Z}_j  + \mathsf U^\top T ((\mathsf U \widetilde{\mathsf  Z}_j) \otimes (\mathsf U \widetilde{\mathsf Z}_j)) +\mathsf U^\top B\mathcal U_j,
$$
where $\widetilde{\mathsf Z}_j$ denotes the $j$-th column of $\widetilde{\mathsf Z}$ and $\mathcal U_j$ denotes the control associated with the dynamical system for the parameter $\mu_j$.
Utilizing the property $ (AB) \otimes (CD) =  (A \otimes C) (B \otimes D)$, we obtain
\begin{equation}
    \label{eq:Zj}
\dot{\widetilde{\mathsf  Z}}_j = A_r \widetilde{\mathsf  Z}_j  + T_r(\widetilde{\mathsf  Z}_j \otimes \widetilde{\mathsf  Z}_j) +B_r \mathcal U_j,
\end{equation}
where\footnote{For the sake of notation, we disregard the dependence of the matrices on the basis $\mathsf U$.} $A_r =  \mathsf U^\top  A \mathsf U \in \mathbb{R}^{r \times r}$, $T_r = \mathsf U^\top T (\mathsf U \otimes \mathsf U) \in \mathbb{R}^{r \times r^2}$ and $B_r = \mathsf U^\top  B \in \mathbb{R}^{r \times m}$.
Equation \eqref{eq:Zj} can be written in the following semilinear form
\begin{equation}
\dot{\widetilde{\mathsf  Z}}_j = \mathcal{A}_r(\widetilde{\mathsf Z}_j) \widetilde{\mathsf Z}_j+ B_r \mathcal U_j,
\label{eq:semi_Z}
\end{equation}
where
$$
(\mathcal{A}_r(H))(i,j) =  A_r(i,j) + \sum_{k=1}^{r} T_r(i,(j-1)r+k) \, H_k ,\quad i,j \in \{1,\ldots, r\},
$$
for $H \in \mathbb{R}^{r}$.
The cost functional \eqref{eq:cost_fun} in the reduced framework reads as 
\begin{equation}
  \min_{\mathcal U_j(\cdot) \in \mathcal{U}_{ad}} J_r(\mathcal U_j(\cdot),\widetilde{\mathsf Z}_j(0)) = 
\int_0^\infty \widetilde{\mathsf Z}_j(s)^\top Q_r \widetilde{\mathsf  Z}_j(s) + \mathcal U_j(s)^\top R \, \mathcal U_j(s) \; ds,
\label{eq:cost_Z}
\end{equation}
with $Q_r = \mathsf U^\top  Q \mathsf U \in \mathbb{R}^{r \times r}$. Given the dynamics \eqref{eq:semi_Z} and the cost functional \eqref{eq:cost_Z}, the associated SDRE reads
\begin{equation}
\mathcal{A}_r(\widetilde{\mathsf Z}_j)^\top P(\widetilde{\mathsf Z}_j) + P(\widetilde{\mathsf Z}_j) \mathcal{A}_r(\widetilde{\mathsf Z}_j) - P(\widetilde{\mathsf Z}_j) F_r P(\widetilde{\mathsf Z}_j) = Q_r,
\label{eq:red_riccati}
\end{equation}
where $P(\widetilde{\mathsf Z}_j) \in \mathbb{R}^{r \times r}$ and $F_r = B_r R^{-1} B_r^\top \in  \mathbb{R}^{r \times r}$.
Upon solving the SDRE, the feedback control can be determined as
$$
\mathcal U_j =  -R^{-1} B_r^{\top} P(\widetilde{\mathsf Z}_j) \widetilde{\mathsf Z}_j.
 $$
 With the feedback representation for the control matrix $\mathcal U$ established, our objective is (ii) to efficiently describe the dynamics in \eqref{eq:UZ_dynamics_feed} \first{in low-rank format}. Let $P_{\widetilde{\mathsf Z}} \in \mathbb{R}^{r \times p}$ be a matrix with columns defined by
 $$
 (P_{\widetilde{\mathsf Z}})(:\,,j) = P(\widetilde{\mathsf Z}_j) \widetilde{\mathsf Z}_j, \quad j=1,\ldots,p.
 $$
 Then,
 $$
\mathsf U^\top (B \mathcal{U}) =  - B_r R^{-1} B_r^\top P_{\widetilde{\mathsf Z}}  = -F_r {P_{\widetilde{\sf Z}}}.
 $$
 
Additionally, we can express
$$
B \mathcal{U} {\mathsf Z}  = -B R^{-1} B_r^{\top} P_{\widetilde{\mathsf Z} } \widetilde{\mathsf Z}^\top = - S (\widetilde{\mathsf Z}).
$$

 Substituting these expressions, we obtain the following DLRA equations:


\begin{equation}
\begin{cases}
    \dot{\mathsf U} = F_{\mathsf U}( \mathsf U,\widetilde{\mathsf Z}) =  (\mathsf I_{N_h} - \mathsf U \mathsf U^{\top}) (A\mathsf U + T ({\mathsf U} \otimes \mathsf U) V(\widetilde{\mathsf Z}) - \second{S(\widetilde{\mathsf Z})\mathsf C^{-1}}), & \\
    \dot{\widetilde{\mathsf Z}} =  F_{\widetilde{\mathsf Z}} (\mathsf U,\widetilde{\mathsf Z}) = A_r \widetilde{\mathsf Z}+ T_r V_{\otimes}(\widetilde{\mathsf Z}) -\second{F_r} {P_{\widetilde{\sf Z}}} , &
\end{cases}
\label{eq:dlra_eq_controlled}
\end{equation}
where
$$
V(\widetilde{\mathsf Z}) = V_{\otimes}(\widetilde{\mathsf Z})\widetilde{\mathsf Z}^\top    \mathsf C^{-1} \in \mathbb{R}^{r^2\times r}
\quad \text{and} \quad
V_{\otimes}(\widetilde{\mathsf Z}) = [\widetilde{\mathsf Z}_1 \otimes \widetilde{\mathsf Z}_1, \ldots,   \widetilde{\mathsf Z}_p \otimes \widetilde{\mathsf Z}_p] \in \mathbb{R}^{r^2 \times p}.
$$
In equations \eqref{eq:dlra_eq_controlled}, we separate the terms that depend on $\mathsf U$ from those that depend on $\widetilde{\mathsf Z}$. This partitioning ensures that the computational cost remains independent of the original dimension $N_h \times p$.

\begin{remark}
\label{remark:celledoni}
Let us point out that the numerical approximation of the system \eqref{eq:dlra_eq_controlled} must keep the orthogonality condition for $\mathsf U$ for all times $t$. This can be achieved in two ways (see, e.g., \cite{lubich06}):
\begin{itemize}
    \item integrating the equation \eqref{eq:dlra_eq_controlled} via a classical method, and then we apply a QR factorization on the new iterate for $\mathsf U$,
    \item applying numerical integration methods able to preserve the orthogonality.
\end{itemize}

We adopt the second approach, drawing on the framework established in \cite{celledoni2002class}, where the authors introduce an intrinsic numerical method designed to preserve the geometric structure of the orthogonal Stiefel manifold, defined as the set of $N_h \times r$ matrices with orthonormal columns. The proposed algorithm achieves a computational complexity of $O(N_h r^2)$, highlighting its efficiency and reliability, and automatically preserves the orthogonality condition for $\mathsf U$.

\end{remark}

\section{Numerical algorithms}
\label{sec:algorithm}
In this section, we present the final algorithms, which incorporate the elements introduced in the previous sections. In Section \ref{sec:POD}, we recall the \first{basic} principles of the POD approach, i.e., the most well-established strategy used to perform model order reduction in this context. In Section \ref{sec:DLRA_cascade}, we introduce the DLRA combined with the \texttt{C-NK} approach to solve the SDRE. In Section \ref{sec:DLRA_control}, we propose the Riccati-based DLRA (R-DLRA) strategies, which guarantee an increased accuracy with respect to the standard DLRA approach counterpart, as we will further discuss in Section \ref{sec:results}. \second{For the sake of clarity, we specify that when we refer to a time discretization $\{t_i\}_{i=1}^{N_T}$, we assume $t_1=0$.}

\subsection{Proper Orthogonal Decomposition}
\label{sec:POD}

This section describes the basic aspects of the POD for nonlinear feedback control. The application of POD within this context is detailed in Algorithm \ref{alg:pod}, where the POD is combined with the \texttt{C-NK} method with cascade information both in the parameters and in time (lines 6, 7, and 10 of Algorithm \ref{alg:pod}). {We remark that Step~10 merely reassigns at time $t_i$ the precomputed value $P_{\mu_1,t_i}$ as the initial guess for the next time step, thereby restarting with the initial guess associated with the first parameter at each time instance.}
The POD approach relies on the construction of fixed basis functions based on an initial computation of the FOM controlled dynamics, called \emph{snapshots}, for a selected subset of parameters $\overline{\mathcal P_h} \subset \mathcal P$ with cardinality $|\overline{\mathcal P_h}| = \overline{p}$ and time instances $\{\overline{t}_i\}_{i=1}^{\overline{N}_T}$. 
More precisely, let us consider a set of solution snapshots denoted as $\second{\{y(\overline{t}_i;\overline{\mu}_j), \; i= 1, \ldots,\overline{N}_T, \; j = 1, \ldots, \overline{p} \}}$, and define the corresponding snapshot matrix as
\begin{equation}
{ S} = [\second{y(\overline{t}_1;\overline{\mu}_1), \ldots, y(\overline{t}_{\overline{N}_T};\overline{\mu}_1), y(\overline{t}_1;\overline{\mu}_2),\ldots, y(\overline{t}_{\overline{N}_T};\overline{\mu}_{\overline{p}})} ] \in \mathbb{R}^{N_h \times N_s},
\label{eq:snap}
\end{equation}
with $N_s = \overline{N}_T \cdot \overline{p}$.
A POD basis of dimension $r \le N_s$ is constructed by performing an orthogonal reduction of the matrix ${ S}$. Specifically, given the Singular Value Decomposition (SVD) of $S$, we obtain
$$
{ S} = { V}{ \Sigma}{ W}^{\top}, \qquad { V}, { W} \in \mathbb{R}^{N_h \times N_s}, { \Sigma} \in \mathbb{R}^{N_s \times N_s}.
$$
The POD basis matrix 
${ V}_r = [{ v}_1 | \ldots |{ v}_r] \in \mathbb{R}^{N_h \times r}$ contains the first $r$ left singular vectors corresponding to the $r$ largest singular values. Given a tolerance $\mathsf{tol}$, the value of $r$ is determined to satisfy the following criterion:
 \begin{equation}
 \label{eq:energy_pod}
     \frac{\sum_{i=1}^r \sigma^2_i}{\sum_{i=1}^{\min(N_h,N_s)} \sigma^2_i} > \mathsf{tol},
      \end{equation}
 where $\{\sigma_i\}_{i=1}^{\min(N_h,N_s)}$ are the singular values of the SVD decomposition.

\begin{algorithm}[H]
\caption{The POD algorithm}
\label{alg:pod}
\begin{algorithmic}[1]
\Statex{\textbf{Input:} $\mathcal{A}_r(\cdot), F, Q,$ tolerance $\mathsf {tol_{\text{\texttt{NK}}}}$, initial conditions $ \mathcal Y_0$, initial guess $P^0_r$, ordered training set $\overline{\mathcal P_h}$, ordered testing set {$\mathcal P_h=\{\mu_i\}_{i=1}^p$}, training time discretization $\{\overline{t}_i\}_{i=1}^{\overline{N}_T}$, time discretization $\{t_i\}_{i=1}^{N_T}$}
\Statex{\textbf{Output:} controlled solutions $\{y^r_{\mu, t_i}\}_{\mu,i}$}
\Statex{}
\State{Compute the snapshot matrix $S$ \second{as in \eqref{eq:snap}}}\label{begin_off}
\State{Perform an SVD of the matrix $S$ and build the POD basis $V_r$} \label{end_off} 
\State{Build \textbf{once and for all} all the reduced matrices \second{appearing in \eqref{eq:sys_POD}}} \label{alg_step_build_pod}
\For{$i = 1, \ldots, N_T-1$}
\For{$\mu \in \mathcal P_h$}
    \State{$P_{\mu, t_i}$ = \texttt{NK}$(\mathcal{A}_r(\cdot), F_r, Q_r, \mathsf {tol_{\text{\texttt{NK}}}}, P_r^0)$}
 \State{$P^0_r = P_{\mu,t_i}$}    
\State{Integrate \eqref{eq:sys_POD} in $(t_i,t_{i+1}]$}
\EndFor
\State{{$P^0_r = P_{\mu_1,t_i}$}}    
\EndFor \label{alg_step_end_pod}
\end{algorithmic}
\end{algorithm}
Using the reduced basis matrix ${ V}_r$, the state vector ${y}_{\mu,t}$ can be approximated as ${y}_{\mu,t} \approx { V}_k y^r_{\mu,t}$, where $y^r_{\mu,t} \in \mathbb{R}^r$ is the solution of the following reduced dynamical system

\begin{equation}
\label{eq:sys_POD}
    \second{\dot{y}^r_{\mu,t}}= A_r y^r_{\mu,t} + T_r (y^r_{\mu,t} \otimes y^r_{\mu,t}) + B_r u(y^r_{\mu,t}),
\end{equation}
where $A_r =  V_r^\top A V_r \in \mathbb{R}^{r \times r}$, $T_r = V_r^\top T(V_r \otimes V_r) \in \mathbb{R}^{r \times r^2}$ and $B_r = V_r^\top  B \in \mathbb{R}^{r \times m}$. We recall that the peculiar quadratic structure \eqref{eq:sys_POD} allows an efficient representation of the reduced system \cite{dolgov2023statisticalproperorthogonaldecomposition,QuarteroniReducedBasisMethods2016}.
{As in the previous section, the reduced system \eqref{eq:sys_POD} can also be expressed in semilinear form as in \eqref{eq:semi_Z}, while the cost functional is similarly reduced to the form in \eqref{eq:cost_Z}, with 
$Q_r = V_r^\top Q V_r$.}
Notice that the reduced basis $V_r$ is fixed: it is precomputed offline and subsequently, the reduced system is assembled and solved for a new parametric instance, in an online phase. 
\subsection{Standard DLRA}
\label{sec:DLRA_cascade}
This section focuses on the DLRA approach for feedback control problems, leveraging the \texttt{C-NK} strategy to solve the SDREs efficiently during the numerical integration process.
The procedure is outlined in Algorithm \ref{alg:dlra}. It begins by performing an SVD of the initial condition matrix $\mathcal Y_0$. 


\begin{algorithm}[H]
\caption{The standard DLRA algorithm}
\label{alg:dlra}
\begin{algorithmic}[1]
\Statex{\textbf{Input:} $\mathcal{A}_r(\cdot), F, Q,$ tolerance $\mathsf {tol_{\text{\texttt{NK}}}}$, initial conditions $ \mathcal Y_0$, initial guess $P^0_r$, ordered set {$\mathcal P_h=\{\mu_i\}_{i=1}^p$}, time discretization $\{t_i\}_{i=1}^{N_T}$}
\Statex{\textbf{Output:} controlled solutions represented as $\{\mathsf U(t_i)\}_{i=1}^{N_T}$ and $\{\mathsf Z(t_i)\}_{i=1}^{N_T}$} 
\Statex{}
\State{SVD($\mathcal Y_0$) $\rightarrow$ $ \mathsf U_0$ and $\mathsf Z_0 = \mathsf U_0^{\top} \mathcal Y_0 $ }
\For{$i = 1, \ldots, N_T-1$}
\label{alg_step_for_dlra}
\State{Build all the reduced matrices \second{appearing in \eqref{eq:sys_POD}}}
\State{$\{P_{\mu,t_i}\}_{\mu \in \mathcal P}$ = \texttt{C-NK}$(\mathcal{A}_r(\cdot), F_r, Q_r, \mathsf {tol_{\text{\texttt{NK}}}}, P_r^0, t_i,\mathcal P_h)$}
\State{Integrate \eqref{eq:dlra_eq_controlled} in $(t_i,t_{i+1}]$}
\State{$P^0_r = P_{\mu_1,t_i}$ }       
\EndFor \label{alg_step_end_dlra}
\end{algorithmic}
\end{algorithm}
Through this process, we obtain the matrix $\mathsf U_0 \in \mathbb R^{N_h \times r}$ and the coefficient matrix $\mathsf Z_0 = \mathcal Y_0^{\top} \mathsf U_0 \in \second{\mathbb R^{p\times r}}$, where $p$ and $r$ are the number of considered parameters and basis functions, respectively.
As for the POD, the number of basis functions $r$ is selected by means of the retained total energy criterion \eqref{eq:energy_pod}.
Therefore, at each time instance, we solve $p$ reduced parametric SDREs \eqref{eq:red_riccati} using the \texttt{C-NK}. The cascade approach is implemented to accelerate the computation of the solution even at the reduced level, as demonstrated in Section \ref{sec:results}. The \texttt{C-NK} considers the information both in the parametric space and in time, except for the very first iteration, i.e., for time $t_1$ and parameter $\mu_1$, the null matrix is employed as an initial guess. Once solved all the SDREs, we assemble the matrices in $F_{\mathsf U}( \mathsf U,\widetilde{\mathsf Z})$ and $ F_{\widetilde{\mathsf Z}}( \mathsf U,\widetilde{\mathsf Z})$ and integrate the system \eqref{eq:dlra_eq_controlled} using a midpoint scheme. Notably, the integration of $\mathsf U$ on the manifold of unitary matrices employs a tangent scheme to ensure geometric consistency, as pointed out in Remark \ref{remark:celledoni}. 
\subsection{Riccati-based DLRA and POD}
\label{sec:DLRA_control}

The DLRA algorithm considers basis functions derived from the SVD of the initial conditions and subsequently evolves them according to the dynamical system \eqref{eq:dlra_eq_controlled}. However, it is noteworthy, that these basis functions are also used for the projection and approximation of the SDREs, even though the initial data does not contain any information related to the Riccati equations. This is not the case with POD, for example, whose construction is based on controlled snapshots based on the full dimensional SDREs. Thus, incorporating additional control-related features in the DLRA approach, such as the Riccati solution or the feedback matrix, could enhance the quality of the approximation. Motivated by the strategies outlined in \cite{schmidt2015basis,schmidt2018reduced}, we propose to employ the information from the Riccati solution ($P$) or the feedback matrix ($K$) to increase the accuracy of the DLRA procedure. We denote these enhanced DLRA approaches with P-DLRA and K-DLRA, respectively. We will refer to these two novel strategies as R-DLRA.

\begin{algorithm}[H]
\caption{The P-DLRA algorithm}
\label{alg:Pdlra}{
\begin{algorithmic}[1]
\Statex{\textbf{Input:} $\mathcal{A}_r(\cdot), F, Q,$ initial conditions $ \mathcal Y_0$, initial guess $P^0_r$, ordered set $\mathcal P_h$, time discretization $\{t_i\}_{i=1}^{N_T}$, number of initial Riccati solutions $s$}
\Statex{\textbf{Output:} controlled solutions $\{\mathsf U(t_i)\}_{i=1}^{N_T}$ and $\{\mathsf Z(t_i)\}_{i=1}^{N_T}$} 
\Statex{}
\State{SVD($\mathcal Y_0$) $\rightarrow$ $\mathsf U_0$ and $ \mathsf U_0^{\top} \mathcal Y_0 = \mathsf  Z_0$ }
\State{Compute $\mathsf P = [P_{\mu_1,0}, \ldots, P_{\mu_s,0}]$ for $s$ parameters}
\State{{$\mathsf  U_P$ = SVD($(I-\mathsf U_0 \mathsf U_0^\top)\second{\mathsf P}$) $\rightarrow \; \; \mathsf U_0 = [\mathsf U_0, \mathsf U_P]$}}
\State{Follow Steps \ref{alg_step_for_dlra}-\ref{alg_step_end_dlra} of Algorithm \ref{alg:dlra}}
\end{algorithmic}
}
\end{algorithm}
Specifically, for the P-DLRA, once built the initial basis $\mathsf U_0$ from $\mathcal Y_0(\mu)$, given a fixed set of parameters $\{\mu_{j}\}_{j=1}^s$ uniformly distributed in $\mathcal{P}$, before the time integration, we compute the corresponding Riccati solutions $\{P_{\mu_{j},0}\}_{j=1}^s$, and collect them in a matrix $\mathsf P = [P_{\mu_1,0}, \ldots, P_{\mu_s,0}] \in \mathbb R^{N_h \times (N_h s)}$. At this point, the basis $\mathsf U_0$ is enriched with a basis $\mathsf U_P$, related to an SVD over $(I-\mathsf U_0 \mathsf U_0^\top)\mathsf P$, to guarantee the orthogonality of the new basis functions with respect to $\mathsf U_0$. Then, the classical DLRA steps are performed in order to solve system \eqref{eq:dlra_eq_controlled}.
An analogous argument can be applied for the K-DLRA, where the initial basis functions in $\mathsf U_0$ are enriched by the basis $\mathsf U_K$ built by an SVD over $(I-\mathsf U_0 \mathsf U_0^\top)\mathsf K$, where 
$\mathsf K = [K_{\mu_1,0}, \ldots, K_{\mu_s,0}] \in \second{\mathbb R^{m\times (N_h s)}}$ collects the feedback matrices related to the Riccati solution in $\mathsf P$, i.e., $K_{\mu_{j},0} = -R^{-1} B^\top P_{\mu_{j},0}$, for $j=1, \cdots, s$.
The P-DLRA and K-DLRA algorithms, which detail this approach, are presented in Algorithms \ref{alg:Pdlra} and \ref{alg:Kdlra}, respectively.
\first{This enrichment procedure is strictly related to the Low-Rank Factor Greedy Algorithm of \cite{schmidt2018reduced} for P-DLRA and to the procedures proposed in \cite{schmidt2015basis} for K-DLRA. Our goal is to exploit such enrichments to increase the accuracy of the DLRA procedure for feedback control, as we will numerically show in Section \ref{sec:results}.}
\begin{algorithm}[H]
\caption{The K-DLRA algorithm}
\label{alg:Kdlra}{{
\begin{algorithmic}[1]
\Statex{\textbf{Input:} $\mathcal{A}_r(\cdot), F, Q,$ initial conditions $ \mathcal Y_0$, initial guess $P^0_r$, ordered set $\mathcal P_h$, time discretization $\{t_i\}_{i=1}^{N_T}$, number of initial Riccati solutions $s$ }
\Statex{\textbf{Output:} controlled solutions $\{\mathsf U(t_i)\}_{i=1}^{N_T}$ and $\{\mathsf Z(t_i)\}_{i=1}^{N_T}$} 
\Statex{}
\State{SVD($\mathcal Y_0$) $\rightarrow$ $\mathsf U_0$ and $\mathsf  U_0^{\top} \mathcal Y_0 = \mathsf Z_0$ }
\State{Compute $\mathsf K = [K_{\mu_1,0}, \ldots, K_{\mu_s,0}]$ for $s$ parameters}
\State{{$\mathsf U_K$ = SVD($(I-\mathsf U_0 \mathsf U_0^\top)\mathsf K$) $\rightarrow \; \; \mathsf U_0 = [\mathsf U_0, \mathsf U_K]$}}
\State{Follow Steps \ref{alg_step_for_dlra}-\ref{alg_step_end_dlra} of Algorithm \ref{alg:dlra}}
\end{algorithmic}
}}
\end{algorithm}
The Riccati-based approach can be employed within the POD algorithm by enriching the basis functions built from the controlled snapshots. For completeness, the P-POD and K-POD are detailed in Algorithms \ref{alg:Ppod} and \ref{alg:Kpod}, respectively. We refer to these two strategies with the common name R-POD.
\begin{algorithm}[H]
\caption{The P-POD algorithm}
\label{alg:Ppod}
\begin{algorithmic}[1]
\Statex{\textbf{Input:} $\mathcal{A}_r(\cdot), F, Q,$ initial conditions $ \mathcal Y_0$, initial guess $P^0_r$, ordered training set $\overline{\mathcal P_h}$, ordered set $\mathcal P_h$, training time discretization $\{\overline{t}_i\}_{i=1}^{\overline{N}_T}$, time discretization $\{t_i\}_{i=1}^{N_T}$, number of initial Riccati solutions $s$ }
\Statex{\textbf{Output:} controlled solutions $\{y_r(t_i,\mu)\}_{i,\mu}$}
\Statex{}
\State{Apply FOM for $t \in \second{ \{\overline{t}_i\}_{i=1}^{\overline{N}_T}}$ and $\mu \in \overline{\mathcal P_h}$ (\texttt{C-NK})}
\State{\second{Follow Steps \ref{begin_off}-\ref{end_off} of Algorithm \ref{alg:pod}} obtaining the basis functions matrix $\mathsf U$}
\State{{\second{Compute $\mathsf P = [P_{\mu_1,0}, \ldots, P_{\mu_s,0}]$ for $s$ parameters}}}
\State{{$\mathsf  U_P$ = SVD($(I-\mathsf U \mathsf U^\top)\mathsf P$) $\rightarrow \; \; V_r = [\mathsf U, \mathsf U_P]$}}
\State{Follow Steps \ref{alg_step_build_pod}-\ref{alg_step_end_pod} of Algorithm \ref{alg:pod} }


\end{algorithmic}
\end{algorithm}

\begin{algorithm}[H]
\caption{The K-POD algorithm}
\label{alg:Kpod}
\begin{algorithmic}[1]
\Statex{\textbf{Input:} $\mathcal{A}_r(\cdot), F, Q,$ initial conditions $ \mathcal Y_0$, initial guess $P^0_r$, ordered training set $\overline{\mathcal P_h}$, ordered set $\mathcal P_h$, training time discretization $\{\overline{t}_i\}_{i=1}^{\overline{N}_T}$, time discretization $\{t_i\}_{i=1}^{N_T}$, number of initial Riccati solutions $s$ }
\Statex{\textbf{Output:} controlled solutions $\{y_r(t_i,\mu)\}_{i,\mu}$}
\Statex{}
\State{Apply FOM for $t \in \second{ \{\overline{t}_i\}_{i=1}^{\overline{N}_T}}$ and $\mu \in \overline{\mathcal P_h}$ (\texttt{C-NK})}
\State{\second{Follow Steps \ref{begin_off}-\ref{end_off}} obtaining the basis functions matrix $\mathsf U$}
\State{\second{Compute $\mathsf K = [K_{\mu_1,0}, \ldots, K_{\mu_s,0}]$ for $s$ parameters}}
\State{{$\mathsf  U_K$ = SVD($(I-\mathsf U \mathsf U^\top)\mathsf K$) $\rightarrow \; \; V_r = [\mathsf U, \mathsf U_K]$}}
\State{Follow Steps \ref{alg_step_build_pod}-\ref{alg_step_end_pod} of Algorithm \ref{alg:pod}}


\end{algorithmic}
\end{algorithm}

\section{Numerical results}
\label{sec:results}
This section presents numerical results comparing the DLRA and R-DLRA methods with standard POD and R-POD reduction techniques. 
\second{The rationale behind the comparison to POD is related to two main aspects: (i) the POD is the most used algorithm when dealing with nonlinear time-dependent problems \cite{QuarteroniReducedBasisMethods2016}, (ii) for the nonlinear time-dependent feedback control problem, finding a reliable error estimator is still an open question and certified reduced basis method has been limitedly employed in this setting.}
In Section \ref{sec:1D}, two experiments are introduced, both centered on an optimal control problem governed by the one-dimensional \both{generalized inviscid} Burgers' equation. \first{The additional advection term in the equation creates a transport-like phenomenon which will highlight the potential of DLRA and R-DLRA in this setting.}
In the first experiment, we demonstrate the advantages of employing \first{the \texttt{C-NK}} strategy in comparison to \texttt{icare} (a direct Matlab solver for Riccati equations) and \texttt{NK} methods across FOM, DLRA, and POD frameworks.
    The second experiment focuses on a comparative analysis of DLRA and R-DLRA against POD and R-POD in terms of computational efficiency (measured by CPU time) and solution accuracy relative to the FOM.
Section \ref{sec:2D} introduces the third experiment, which extends the comparison to a two-dimensional optimal control problem constrained by the \both{generalized inviscid} Burgers' equation. Here, the performance of DLRA and R-DLRA is evaluated against POD and R-POD, considering both computational efficiency and accuracy with respect to the FOM solution.

\subsection{1D \both{generalized inviscid} Burgers' equation} 
We define the computational domain $\Omega=[-5,30]$ and the computational time interval $I=[0,T]$, with $T=1$. Throughout the following discussion, $\chi_{\omega}(\cdot)$ denotes the indicator function over the subdomain $\omega$. The interval $\omega_c$ specifies the region where the control action is active, while $\omega_o$ indicates the interval where we observe the dynamics. At the continuous level, the optimal control problem is formulated as follows:
\begin{equation}
\label{eq:func_bur}
    J(u,y_0) = \int_0^{+\infty} \int_{\omega_0} | y(t,x)|^2 \, dx + 0.1 \int_{\omega_c} | u(t,x) |^2 \, dx \, dt,
\end{equation}
subject to 
\begin{equation}
    \frac{d}{dt}y(t,x) + y(t,x) \partial_{x} y(t,x) + 20 \partial_{x} y(t,x) = u(t,x)\chi_{\omega_c}(x), 
\end{equation}
for a given parametric initial condition $y_{0}(x;\mu) = \mu_1 e^{-\mu_2 x^2}$ with $ \mu = (\mu_1, \mu_2) \in \mathcal P = [0.1,0.5]\times [0.2,1.5]$. Homogeneous Neumann boundary conditions are imposed for $x=-5$ and $x=30$. For the numerical discretization, a Finite Difference scheme in space is applied, resulting in a FOM of dimension $N_h = 100$. Time integration is performed using a midpoint quadrature rule with a time step of $\Delta t = 10^{-3}$. In this study, 
the control acts on the subdomain $\omega_c = [-1,1] \cup [5,7]$, while the observation is taken over $\omega_o = [2,4] \cup [8,10]$. The POD and R-POD models are built with a $4 \times 4$ parametric grid over $\mathcal P$, with 100 snapshots in time for each parameter, resulting in a total of $N_{s} =1600$ snapshots. The DLRA and R-DLRA approaches compute the solution at all time instances using a $20 \times 20$ parametric grid over $\mathcal P$, corresponding to $p = 400$ parameters. 
\first{In the \texttt{C-NK} setting, the parameters are ordered according to the following rule: we first vary the first parameter while keeping the second parameter fixed. Specifically, starting with $j =1$, we consider the pairs $(\mu^j_1,\mu^1_2), \ldots, (\mu^j_1,\mu^{20}_2)$. Subsequently, we increment $j$ to $j+1$ and repeat this process until $j=20$.
}
This parametric grid is also used as a testing set for the POD method. \first{We remark that exploiting the results of Section \ref{sec:sensitivity} to compute $\Delta \mu$ and $\Delta t$ is impractical since the bound should be computed for every parameter and time instance. 
The results of Section \ref{sec:sensitivity} are used as a \emph{starting point} to choose a $\overline{\Delta \mu}$ and $\overline{\Delta t}$ for the first parameter and initial time, respectively. The values $\Delta \mu$ and $\Delta t$ we used in the numerical simulations were smaller than $\overline{\Delta \mu}=1.15$ and $\overline{\Delta t}=0.08$. A posteriori, no convergence issues were observed for the \texttt{C-NK}. The same procedure was exploited in Section \ref{sec:2D}}. The reduced dynamics for both DLRA and POD strategies are integrated consistently with the FOM. The performances are tested in terms of CPU times and $L^2$ absolute error in the \first{state dynamics} and the control for  $t \in I$. The state and control errors are denoted by 
\begin{equation*}
E = E(t) =  \frac{1}{|\mathcal P_h|}{\sum_{\mu \in \mathcal{P}_h}\norm{y_{\mu,t} - y^r_{\mu,t}}_{L^2(\Omega)}} \quad \text{and} \quad 
E_c = E_c(t) =  \frac{1}{|\mathcal P_h|}{\sum_{\mu \in \mathcal{P}_h}\norm{u_{\mu,t} - u^r_{\mu,t}}_{L^2(\omega_c)}},    
\end{equation*}
    respectively, where we recall that $y^r_{\mu,t}$ represent the reduced (i.e., DLRA or POD) state solution at time $t$ for a given $\mu$, while $u^r_{\mu,t}$ represent the reduced (i.e., DLRA or POD) control solution at time $t$ for a given $\mu$, to be compared to the FOM solutions. 
\label{sec:1D}

\medskip

\textbf{Experiment 1: What is the best option between \texttt{icare}, \texttt{NK}, and \texttt{C-NK}?} 

\medskip

In this test, we investigate the performances of the Matlab function \texttt{icare}, \texttt{NK} scheme (Algorithm \ref{alg:newton}), and \texttt{C-NK} method (Algorithm \ref{alg:cnk}) to solve the SDRE related to the optimal control problem of Section \ref{sec:1D} both at the FOM and ROM level. For DLRA, we need to apply an SVD over the parametric initial condition $\mathcal Y_0(\mu) \in \mathbb R^{N_h \times p}$. We choose the basis functions by an energy criterion \eqref{eq:energy_pod} over the parametric initial condition, with a tolerance $\mathsf{tol} = 0.9999$.
DLRA meets the criterion for $r=4$.
For the POD, we retain the same number of basis functions for an SVD over the $N_{s}$ snapshots, to be comparable to the DLRA setting.

The indicators taken into consideration for the comparison are:
\begin{itemize}
\item the {global CPU} time, i.e., the time needed to perform the simulations for the $p=400$ parameters. For the POD, the global CPU also considers the offline phase \cite{Pagliantini2021409}. This is a natural choice: the DLRA procedure is \emph{on the fly} and does not need any training compared to the POD approach. Thus, the global CPU should consider the offline time, which is a mandatory step for the POD;
\item the CPU of \texttt{NK}, i.e., the CPU time, averaged in time and in the parametric space, needed for the \texttt{NK}-based strategies to converge;
\item the average number of iterations needed to the \texttt{NK}-based strategies to converge;
\item the average state and control error. 
\end{itemize} 
The results are summarized in Table \ref{tab:test1}. The initial guess for \texttt{NK} is consistently the trivial one (the null matrix). In contrast, \first{the \texttt{C-NK}} utilizes the trivial guess solely at the first time instance and for the first parameter. Subsequently, it employs cascade information guided by the sensitivity analysis discussed in Section \ref{sec:sensitivity}. \\
We observe that \first{the \texttt{C-NK}} obtains the lowest global \first{CPU time} for \emph{all} the strategies. While this was quite expected at the FOM level (and consequently for the POD during the offline stage), this conclusion was not obvious for DLRA. Indeed, DLRA \emph{only solves reduced Riccati} and direct methods, such as \texttt{icare}, can be highly competitive. However, from the table, we note that \texttt{C-NK} is two times faster than the direct method and slightly faster than the classical \texttt{NK}. Indeed, using the cascade information reduces the time of convergence of the method and the number of iterations needed for the three algorithms. In addition, for the reduced strategies, the accuracy achieved with the iterative solver is the same of the one obtained with the \texttt{icare} direct solver \first{up to the eleventh significant digit}. Based on these results, we conclude that \texttt{C-NK} is the most suitable choice for the FOM, DLRA, and POD-based strategies in the subsequent experiments.
    \begin{table}[H]
    \begin{center}
        
\caption{\centering Experiment 1. CPU and accuracy comparison for \texttt{icare, NK} and \texttt{C-NK} between FOM, DLRA and POD ($r=4$).}
\vspace{3mm}
\label{tab:test1}
\resizebox{\textwidth}{!}{
\centering
\begin{tabular}{c|ccc|ccc|ccc}
           & \multicolumn{3}{c|}{\texttt{icare}}                       & \multicolumn{3}{c|}{\texttt{NK}} & \multicolumn{3}{c}{\texttt{C-NK}} \\[1.5mm]
         & FOM       & DLRA  & POD   & FOM     & DLRA  & POD     & FOM     & DLRA    & POD    \\ \hline
Global CPU      & $1.45\cdot 10^4$ & $241$ & $843$ & $1.85\cdot 10^4$ & $138$ & $748$ & $8.72\cdot 10^3$ & $122$   & $496$  \\ \hline
CPU of \texttt{NK}      & - & -  & - & $3.55\cdot 10^{-2}$ & $2.34\cdot 10^{-4}$ & $3.46\cdot 10^{-4}$ & $1.92\cdot 10^{-2}$ & $2.11\cdot 10^{-4}$   & $2.51\cdot 10^{-4}$  \\ \hline

Iterations & - & - & - & $4$      & $3.61$     & $5$     & $1.29$     & $1.34$    & $1.81$    \\ \hline
Average $E$ & -         & $1.10\cdot 10^{-3}$    & $5.10\cdot 10^{-3}$   & -       & $1.10\cdot 10^{-3}$   &  $5.10\cdot 10^{-3}$     & -       & $1.10\cdot 10^{-3}$ & $5.10\cdot 10^{-3}$ \\ \hline
Average $E_c$ & -         & $3.91\cdot 10^{-4}$ & $3.10\cdot 10^{-3}$   & -       & $3.91\cdot 10^{-4}$ & $3.10\cdot 10^{-3}$      & -       & $3.91\cdot 10^{-4}$ & $3.10\cdot 10^{-3}$ \\ \hline

\end{tabular}}
\end{center}
\end{table}

\medskip

\textbf{Experiment 2: What is the best option between DLRA, R-DLRA, POD, and R-POD?}

\medskip

\begin{figure}[H]	
\begin{center}    
	\includegraphics[scale=0.5]{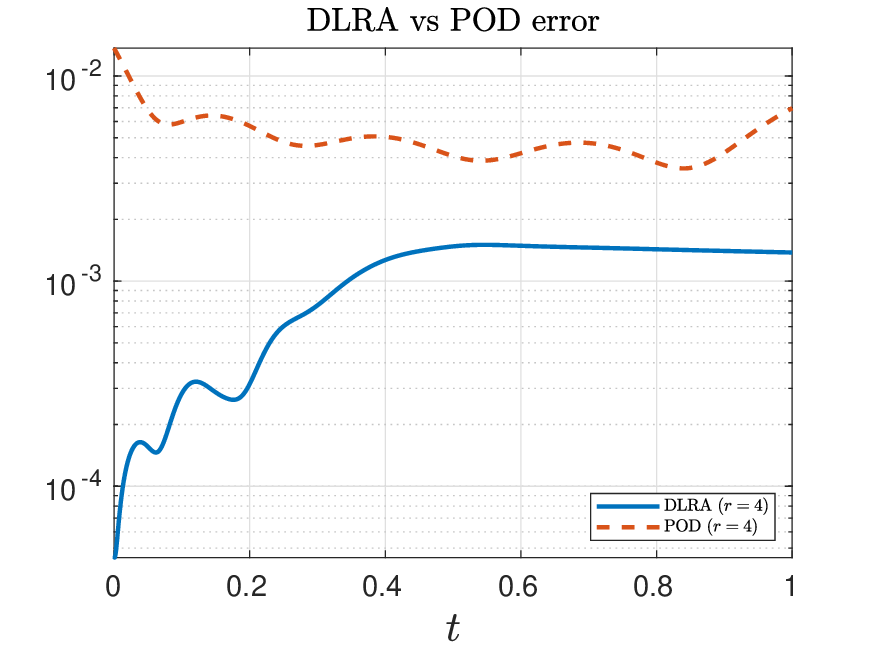}	
    \includegraphics[scale=0.5]{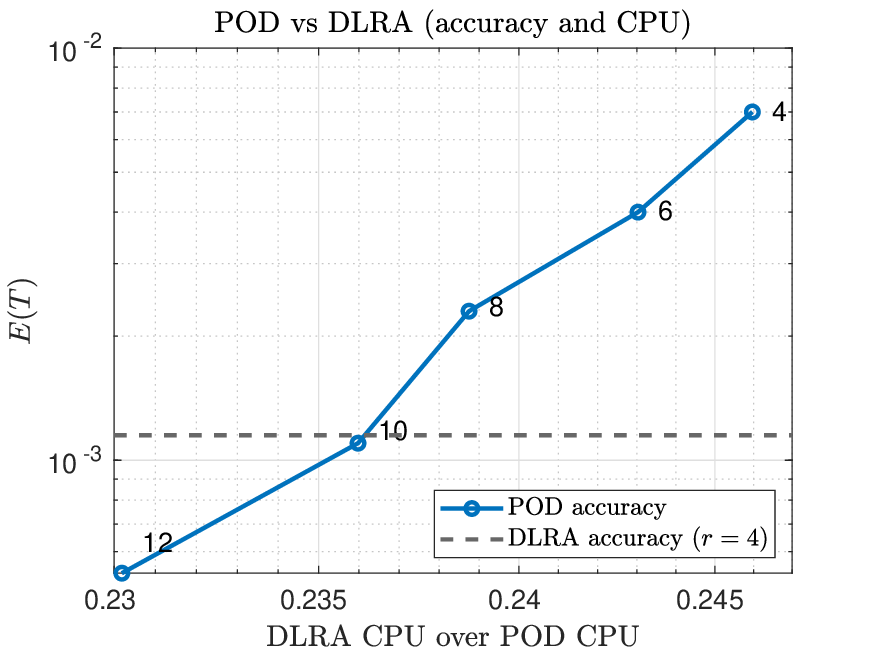}
\caption{Experiment 2. \emph{Left}: Evolution of $E(t)$ for DLRA and POD for $r=4$. \emph{Right}: POD final accuracy plot with CPU times for increasing values of $r$ to be compared to the final accuracy of the DLRA.}
\label{pics:1d_DLRAvsPOD}
\end{center}

\end{figure}

In this test case, we investigate the role of the Riccati information in terms of CPU and the accuracy of the DLRA and the POD approaches compared to the FOM solution. The SDREs are always solved using the \texttt{C-NK} algorithm, which proved to be the best approach both at the FOM and ROM settings in Experiment 1. Let us start our analysis from DLRA and POD. We are working in the numerical setting of the previous test case, with $r=4$ basis functions. In Figure \ref{pics:1d_DLRAvsPOD}, we compare the DLRA and POD approaches in terms of the evolution of the state error $E(t)$ and the CPU times. 
In the left panel, the evolution of $E$ in time is represented. For the same $r$, the DLRA approach provides a more accurate recovery of the FOM dynamics over time compared to the POD method. In the right panel, we fix $r=4$ for the DLRA, leading to a final accuracy of $E(T) = 1.4\cdot 10^{-3}$ and compare it to the accuracy of the POD method for increasing values of $r$. We also report the ratio of the CPU time for DLRA to the CPU time for POD. To achieve the same accuracy as DLRA, the POD method requires $r=10$ basis functions. However, DLRA proves to be significantly more computationally efficient, being approximately four times faster. The same conclusions can be drawn by analyzing the first column of Table \ref{tab:global_tab}: DLRA emerges as the most efficient approach, demonstrating a lower CPU time and a significant \emph{speedup} compared to the FOM. Here, \emph{speedup} is defined as the ratio of the CPU time required by the FOM to solve the dynamics for $p=400$ parameters to the CPU time required by the reduced-order strategy to accomplish the same task.

\begin{figure}[H]	
\begin{center}

\includegraphics[scale=0.5]{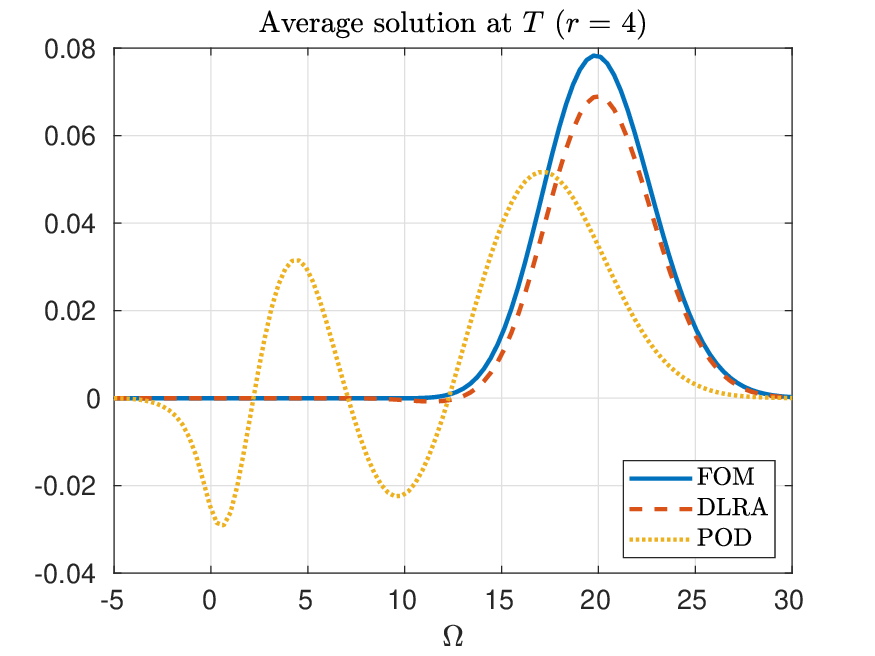}
	\includegraphics[scale=0.5]{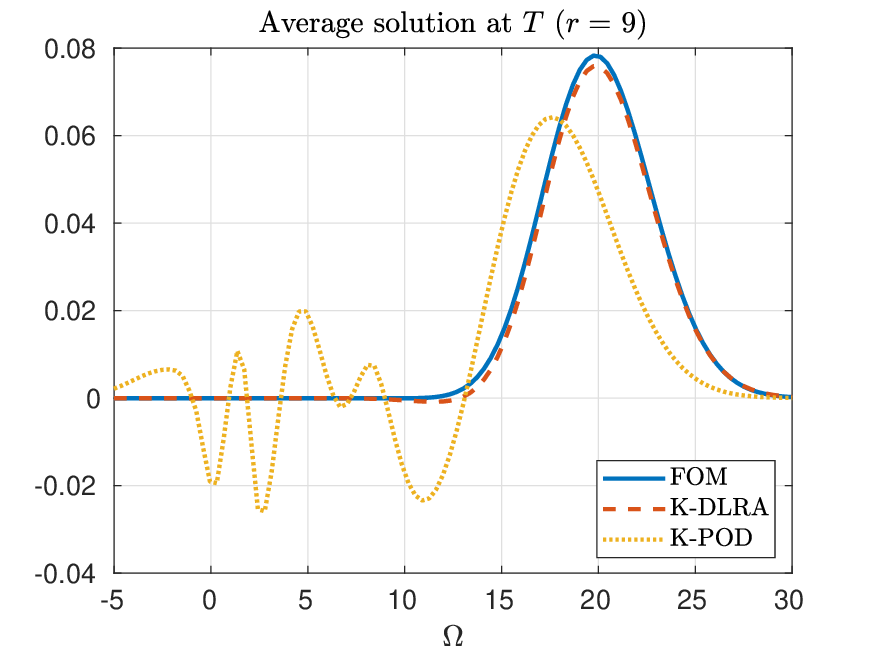}	
	\includegraphics[scale=0.5]{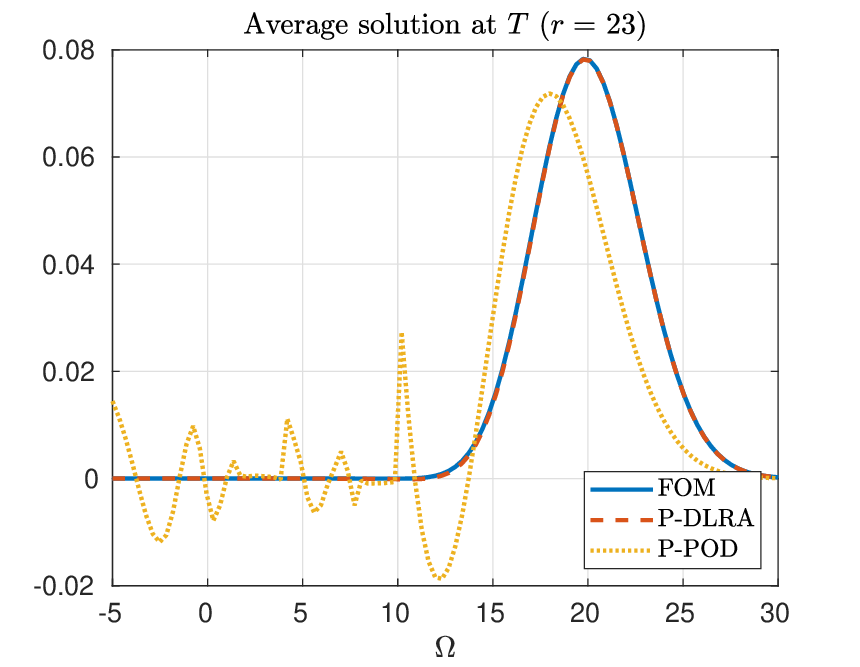}	
        
\caption{Experiment 2. Solutions averaged in $\mathcal P_h$ at the final time $T$. \emph{\first{Top left}}: FOM, DLRA, and POD for $r=4$. \emph{\first{Top right}}: FOM, K-DLRA, and K-POD for $r=9 \; (r_2 = 5)$. \emph{Bottom}: FOM, P-DLRA, and P-POD for $r=23 \; (r_2 = 19)$.}
\label{pics:solutions1D}
\end{center}

\end{figure}

The superior accuracy of the DLRA approach is illustrated in Figure \ref{pics:solutions1D}. The top plot compares the FOM, DLRA, and POD solutions averaged over $\mathcal P_h$ at the final time $T$. While the POD method fails to accurately reconstruct the solution, the DLRA approach captures the dynamics of the FOM more effectively. However, it is noted that the magnitude of the solution is not fully recovered by the DLRA. 

\begin{figure}[H]	
\centering
\includegraphics[scale=0.5]{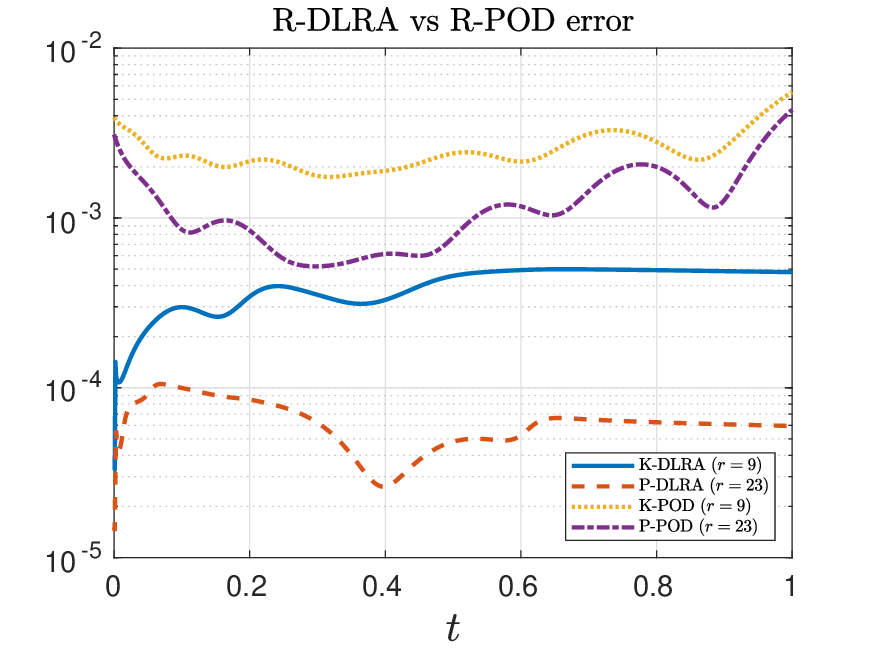}\\
	\subfloat[\centering \first{Comparison varying $r_2$}]{\includegraphics[scale=0.5]{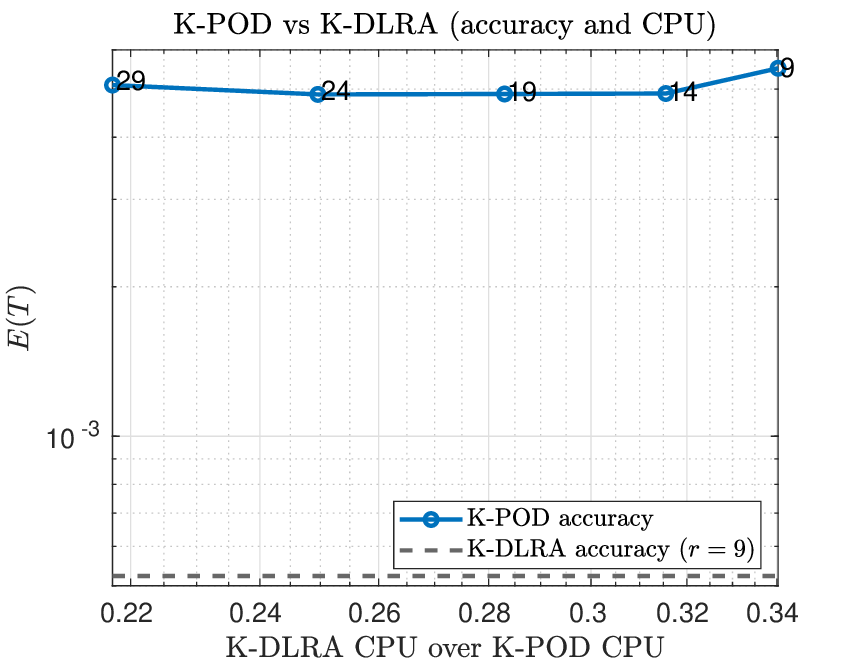} \label{a1}}	
    \subfloat[\centering \first{Comparison varying $r_1$}]{
    \includegraphics[scale=0.5]{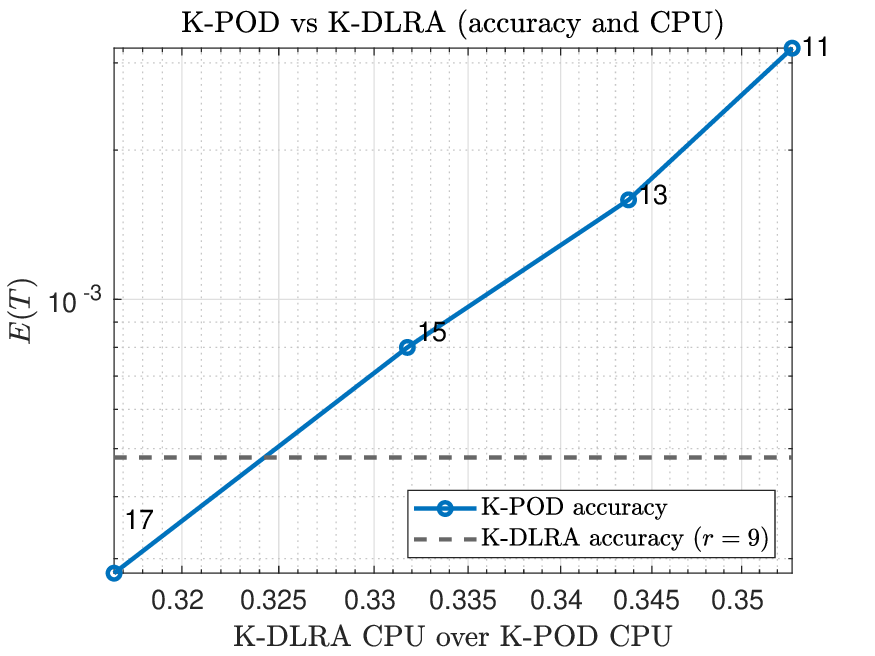} \label{b1}}	\\
    \subfloat[\centering \first{Comparison varying $r_2$}]{
	\includegraphics[scale=.5]{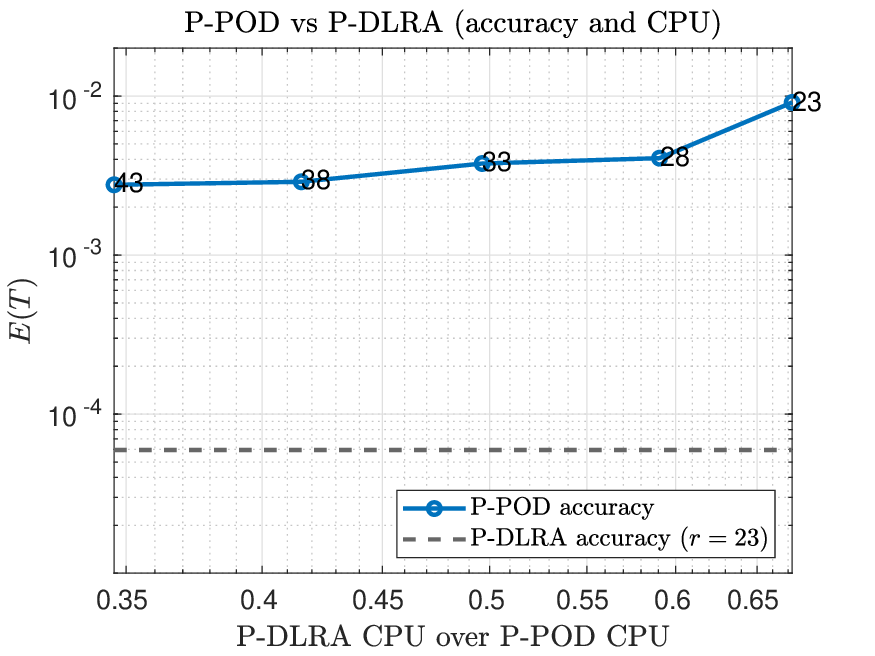}\label{c1}}
    \subfloat[\centering \first{Comparison varying $r_1$}]{
	\includegraphics[scale=0.5]{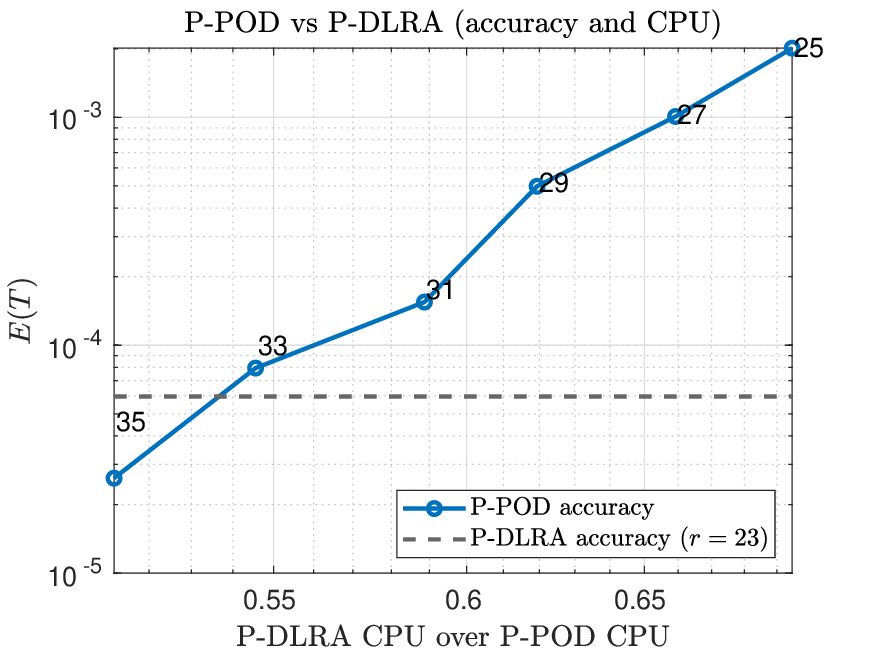}\label{d1}}
	
\caption{Experiment 2. \emph{Top}: Evolution of $E(t)$ for K-POD, K-DLRA, P-POD, and P-DLRA with $r=9$ and $r=23$, respectively. \emph{Middle left}: K-POD final accuracy with CPU times for increasing values of $r_2$ to be compared to the final accuracy of the K-DLRA \protect\subref{a1}. \emph{Bottom left}: analogous representation for P-POD \protect\subref{c1}. \emph{Middle right}: K-POD final accuracy with CPU times for increasing values of $r_1$ to be compared to the final accuracy of the K-DLRA \protect\subref{b1}. \emph{Bottom right}: analogous representation for P-POD \protect\subref{d1}. }
\label{pics:KPerrs}
\end{figure}

We overcome this issue thanks to the R-DLRA approaches. When dealing with Riccati-based strategies, we work with $r=r_1 + r_2$ basis functions, where $r_1$ is the number of basis functions retained by the SVD of the initial state, i.e., $r_1=4$, while $r_2$ is the number of additional bases related to the Riccati information computed {for the first parametric and time instace: only one Riccati was sufficient to achieve a drastic improvement in terms of accuracy.}
 In the case of K-DLRA, $r_2=5$ retains the 0.9999 of the energy, while for P-KLRA, to retain the same amount of energy, $r_2=19$ is needed, leading to a global basis number of $r=9$ and $r=23$ for the two approaches, respectively. From the bottom row of Figure \ref{pics:solutions1D}, we see that adding Riccati bases is undoubtedly beneficial for the DLRA. By employing this approach, we closely track the dynamics of the system, achieving a visual overlap between the FOM and the P-DLRA profiles. This is not the case for K-POD and P-POD: the Riccati information yields spurious oscillations in the solution. \\   
In Figure \ref{pics:KPerrs}, the top panel shows that P-POD improves the accuracy of the solution compared to K-POD and standard POD (as indicated by the red dashed line in the left plot of Figure \ref{pics:1d_DLRAvsPOD}). However, P-POD provides a less accurate representation of the FOM solution when compared to P-DLRA (it is even less accurate than K-DLRA). 
In the same figure, in the middle and bottom plots, we propose the performance analysis in terms of accuracy and CPU times. In the left plots, we fix $r$ for the K-DLRA and P-DLRA, while we fix $r_1=4$ and we increase $r_2$ for the R-POD. It is evident that adding more Riccati-based \first{basis functions} does not improve the accuracy of the POD method. The only way to reduce the error is to fix $r_2$ ($r_2=5$ and $r_2=19$ for K-POD and P-POD, respectively) and increase $r_1$. Once again, while R-POD strategies can achieve the same accuracy as their R-DLRA counterparts, they require significantly more computational time. Specifically, K-DLRA is nearly three times faster than K-POD, and P-DLRA is approximately two times faster than P-POD.

Table \ref{tab:global_tab} summarizes the CPU time for the Riccati-based strategies. Once again, K-DLRA and P-DLRA outperform their respective POD counterparts. While P-DLRA demonstrates superior accuracy in capturing the model dynamics, it does so at the cost of increased computational time, as expected. Nevertheless, DLRA, K-DLRA, and P-DLRA consistently require less computational time compared to their corresponding POD-based approaches.

\second{We remark that the CPU times for R-DLRA and R-POD take into account the time needed to construct the additional matrices ($P$ or $K$). DLRA presents a short offline phase, which builds the initial basis $\mathsf U_0$ and $\mathsf Z_0$ from the initial condition $\mathcal Y_0$. In this setting, the computation of the additional Riccati information highly influences the CPU time with respect to standard DLRA. This is not the case for POD, which presents a much heavier offline phase based on the controlled trajectory compression. In our numerical results, this leads to a less dramatic difference between the POD and the R-POD CPU times. Moreover, we highlight that DLRA has a heavier online phase with respect to POD, since the basis functions must be updated at each time instance.}
\begin{table}[H]
\caption{\centering Experiment 2. CPU comparison for DLRA, POD, K-DLRA, K-POD, P-DLRA and P-POD.}
\label{tab:global_tab}
\centering
\begin{tabular}{c|cc|cc|cc}
            & \multicolumn{2}{c|}{$r=4$} & \multicolumn{2}{c|}{$r=9$} & \multicolumn{2}{c}{$r=23$} \\
            & DLRA         & POD         & K-DLRA       & K-POD       & P-DLRA       & P-POD       \\ \hline
global CPU  & $122$          & $496$         &   $173$        & $509$         &    $458$       & $680$         \\ \hline
speedup     & $71$           & $17$          &     $50$       &    $17$         &     $19$     &    $12$         \\ \hline
\end{tabular}
\end{table}
\begin{figure}[H]	
\centering
	\includegraphics[scale=0.5]{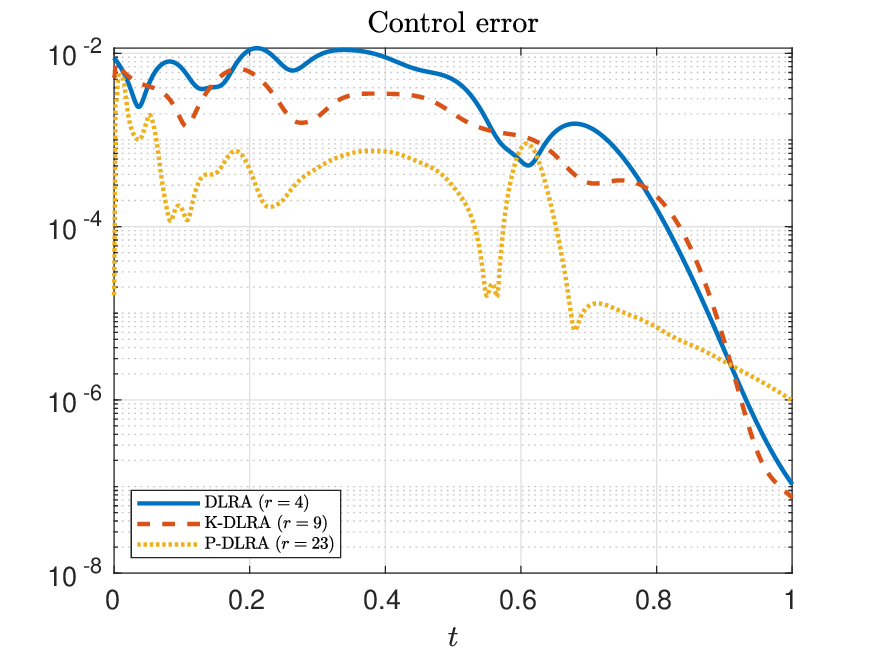}	
	\includegraphics[scale=0.5]{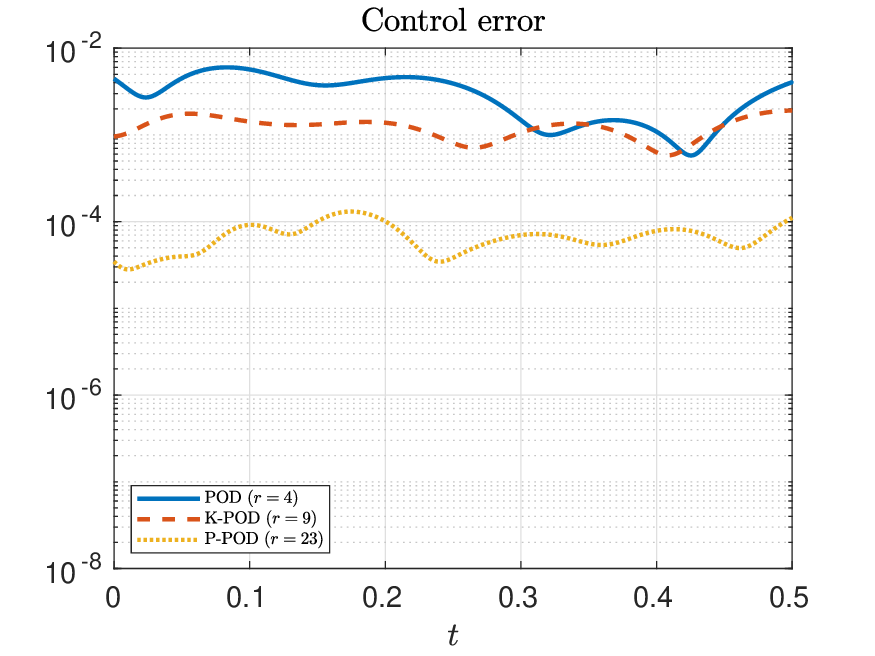}
\caption{Experiment 2. Evolution of $E_c(t)$. \emph{Left plot}: DLRA, K-DRLA, and P-DLRA. \emph{Right}: analogous plot for POD, K-POD, and P-POD. }
\label{pics:err_control}
\end{figure}
For the sake of completeness, Figure \ref{pics:err_control} presents the evolution of the control error \first{$E_c(t)$. In} the left panel, we depict the control error for DLRA, K-DLRA, and P-DLRA, while the right panel presents a similar plot for POD, K-POD, and P-POD. Initially, DLRA and K-DLRA show similar performance to POD and K-POD at the beginning of the simulation. However, after $t=0.7$, the control error for DLRA and K-DLRA decreases significantly, outperforming the performance of POD and K-POD. A similar trend is observed for P-DLRA: although P-POD initially benefits from the Riccati information, by $t=0.65$,  P-DLRA becomes more accurate than P-POD. On average, in terms of control accuracy, POD and K-POD perform similarly to DLRA and K-DLRA. P-POD. However, P-POD performs better than P-DLRA in terms of control accuracy, on average. Nonetheless, it is important to note that POD, K-POD, and even P-POD fail to adequately represent the dynamics, as shown in Figure \ref{pics:solutions1D}, leading to a degradation in control accuracy in time. 

\begin{remark}[The role of the control]
\label{rem:control}
    In Figure \ref{pics:err_control}, we notice that the DLRA strategies, at the beginning of the simulation, struggle to represent the control action. This might seem unexpected, but we recall that the DLRA process is not an \emph{explore-compress} algorithm over the controlled trajectory, but it only relies on the initial state, with scarce (e.g., K-DLRA and P-DLRA) or no control information (DLRA).
    In contrast, the POD strategies have significantly more information regarding the control action, as the basis functions are constructed using the \emph{controlled dynamics}. This additional control-related information can aid in more accurately recovering the control action at the beginning of the simulation.  
    \end{remark}

\begin{figure}[H]	
\centering
\includegraphics[scale=0.5]{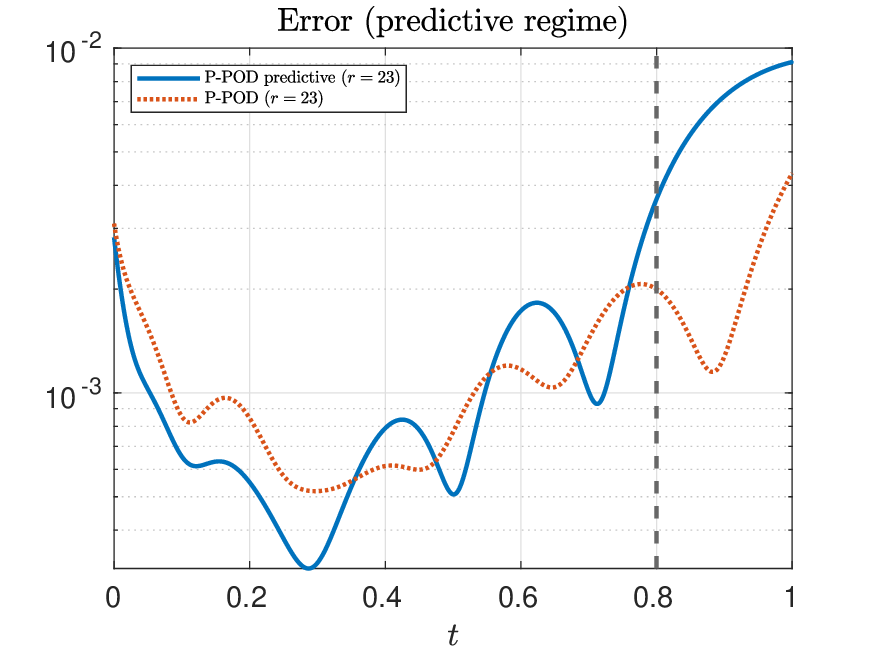}
\includegraphics[scale=0.5]{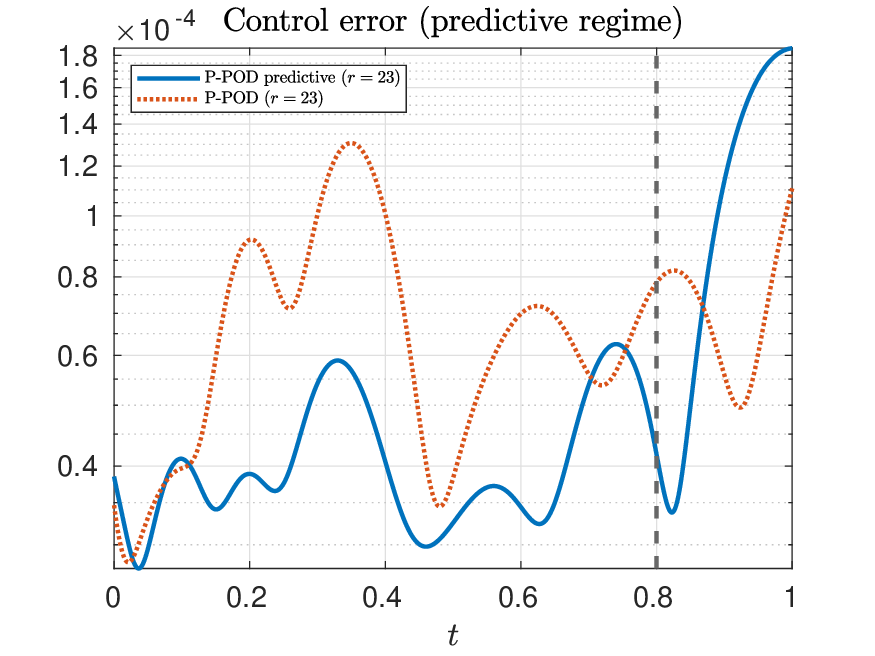}
\caption{Experiment 2, predictive regime for P-POD. \emph{Left plot}: evolution of $E(t)$ for P-POD in the predictive regime and P-POD. \emph{Right plot}: analogous plot for $E_c(t)$.}
\label{pics:predictive}
\end{figure}

\begin{remark}[The actual advantage of DLRA]
Thus far, we have focused on the reconstructive regime of POD. However, it is important to emphasize that DLRA and its Riccati-based variants operate as on-the-fly processes, eliminating the need for a training phase. This represents a significant advantage over POD strategies, which may fail in a \emph{predictive} setting. For completeness, we conduct a simulation using P-POD (the most effective POD-based reduction method), where the snapshots are taken from the same dynamics as the previous test case but limited to $t=0.8$, with the number of snapshots fixed at $N_{s}=1600$. As expected, the accuracy of P-POD significantly deteriorates towards the end of the simulation, where no prior information is available. Indeed, from the left plot of Figure \ref{pics:predictive} we observe that in the predictive regime, P-POD performs worse than its reconstructive counterpart. By comparison with Figure \ref{pics:KPerrs}, we can also notice that the error for the P-DLRA does not increase in time. This behavior is unsurprising, given that DLRA does not rely on training or require prior system information. The same conclusions are drawn from the right plot of Figure \ref{pics:predictive}, which depicts the control error for P-POD in both the predictive and reconstructive regimes. Once again, we note that by the end of the simulation, P-DLRA outperforms P-POD in representing the control even in the reconstructive regime, see Figure \ref{pics:err_control}. For these reasons, in the next test, we only focus on the reconstructive regime, the only setting where the POD strategies might be competitive compared to the DLRA approaches.         

\end{remark}

    \first{
\begin{remark}[Evaluating SDRE Optimality]
{We now address} the question of how closely the SDRE controller approximates the optimal control. 
In this {specific} scenario, an exact computation of the optimal control is not feasible because it would require solving the infinite‑dimensional Hamilton–Jacobi–Bellman equation associated with the nonlinear system, which suffers from the curse of dimensionality and admits no closed‑form solution except in trivial settings. As an alternative, the Pontryagin Maximum Principle (PMP) is employed to compute an open-loop control by solving the associated state-adjoint system over a sufficiently long time horizon. The resulting control inputs and state trajectories are then compared with those generated by the SDRE controller.
We compare the methods in the time interval $[0,1]$, while the PMP system will be derived in the larger time interval $[0,2]$. The final time $t=2$ has been selected such that further increases in the time horizon do not produce significant changes in the resulting solutions. The corresponding PMP system is solved via a gradient descent approach, wherein the state-adjoint system is iteratively integrated starting from an initial guess for the control input. In this study, the initial control is set to zero. The gradient-based optimization is performed using the MATLAB function \texttt{fminunc}. Additional methodological details can be found in \cite{azmi2020optimal} and the references therein. 
For the comparison {with respect to the Riccati-based feedback controlled FOM solution}, we utilize error indicators analogous to those previously defined, denoted here as $\tilde{E}$ and $\tilde{E}_c$, given by:
\begin{equation*}
\tilde{E} =  \frac{1}{|\mathcal P_h|}{\sum_{\mu \in \mathcal{P}_h}\norm{y_{\mu,t} - y^{PMP}_{\mu,t}}_{L^2(\Omega)}} \quad \text{and} \quad 
 \tilde{E}_c =  \frac{1}{|\mathcal P_h|}{\sum_{\mu \in \mathcal{P}_h}\norm{u_{\mu,t} - u^{PMP}_{\mu,t}}_{L^2(\omega_c)}},    
\end{equation*}
where $y^{PMP}_{\mu,t}$ and $u^{PMP}_{\mu,t}$ denote, respectively, the trajectory and control computed via the PMP approach at time $t$ for a given parameter $\mu$. The computed results yield an average $\tilde{E}$ of $7.88 \cdot 10^{-5}$ and an average $\tilde{E}_c$ of $1.98 \cdot 10^{-4}$. From these errors, we deduce that the feedback control approximation is reliable in terms of optimality, and the reduced approximation errors of Table \ref{tab:test1} dominate with respect to the optimality error we are committing.
\end{remark}
}

\subsection{2D \both{generalized inviscid} Burgers' equation}
\label{sec:2D}

We define the computational domain $\Omega=[-10,10] \times [-10,10]$ and the time interval $I=[0,T]$, with $T=0.5$. The optimal control problem we solve is:
\begin{equation}
\label{eq:func_bur}
    J(u,y_0) = \int_0^{+\infty} \int_{\omega_o}| y(t,x)|^2 \, dx + 0.1 \int_{\omega_c}| u(t,x) |^2 \, dx \, dt,
\end{equation}
subject to 
\begin{equation}
    \frac{d}{dt}y(t,x) + y(t,x) \nabla   \cdot y (t,x) + 20 \nabla \cdot y(t,x) = u(t,x)\chi_{\omega_c}(x), 
\end{equation}
where $\omega_c = [-5,0] \times [-5,0]$ is the control domain, $\omega_o = [0,5] \times [0,5]$ is the observation domain, and $x = (x_1, x_2)$ is the spatial variable.

The parametric initial condition is $y_{0}(x;\mu) = \mu_1 e^{-\mu_2 [(x_1+2)^2 + (x_2+2)^2]}$ for $ \mu = (\mu_1, \mu_2) \in \mathcal P = [0.01,0.05]\times [0.1,0.3]$. Homogeneous Neumann boundary conditions are imposed for the boundary $\partial \Omega$. For the space discretization, we apply a Finite Difference scheme. The dimension of the FOM is $N_h = 400$. In time, we perform a midpoint integration scheme with $\Delta t = 10^{-3}$. As in the previous case, the POD and R-POD models are built with a $4 \times 4$ parametric grid over $\mathcal P$, with 100 snapshots in time for each parameter, resulting in a total of $N_{s} =1600$ snapshots. The DLRA and R-DLRA approaches compute the solution at all time instances using a $20 \times 20$ parametric grid over $\mathcal P$, corresponding to $p = 400$ parameters. This set of parameters is the testing parametric set for the POD methods. The order of the parameters for the \texttt{C-NK} strategy follows the one proposed in Section \ref{sec:1D}.

The performances are tested in terms of CPU times, $E(t)$ and $E_c(t)$. The basis functions for DLRA and R-DLRA have been selected using the same energy criterion as in Experiment 1 and Experiment 2. The criterion leads to $r=3$ for DLRA and $r= r_1 + r_2 = 11$ and $r=r_1+r_2=41$, with $r_1=3$ and $r_2=9$ and $r_2=38$, for K-DLRA and P-DLRA, respectively. We recall that $r_1$ is the number of basis \first{functions} from the SVD of the initial state, while $r_2$ is the number of basis functions related to the additional Riccati information, computed in the same parametric set of the standard DLRA approach.

\begin{figure}[H]	
\centering
\includegraphics[scale=0.33]{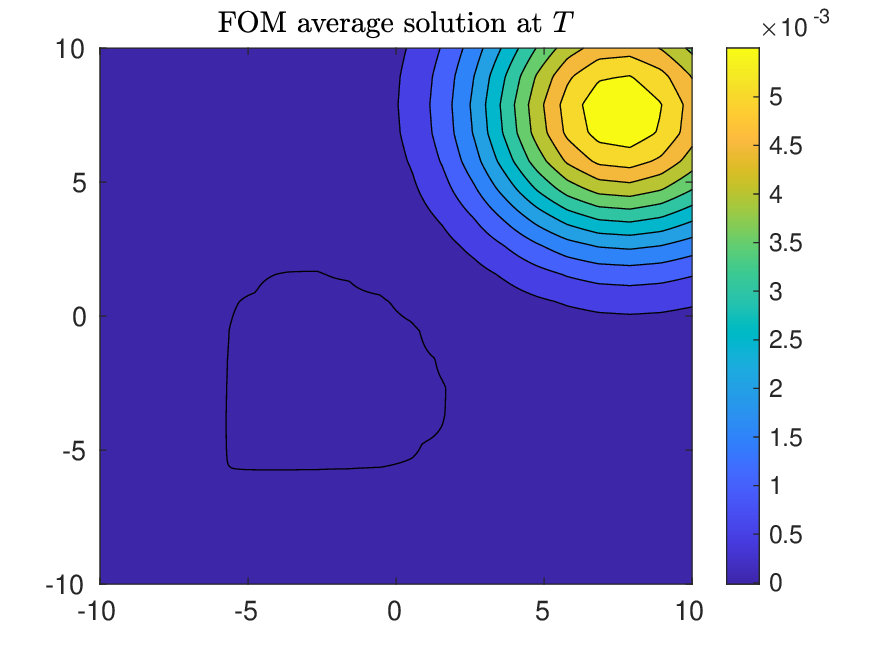}\\
	\includegraphics[scale=0.33]{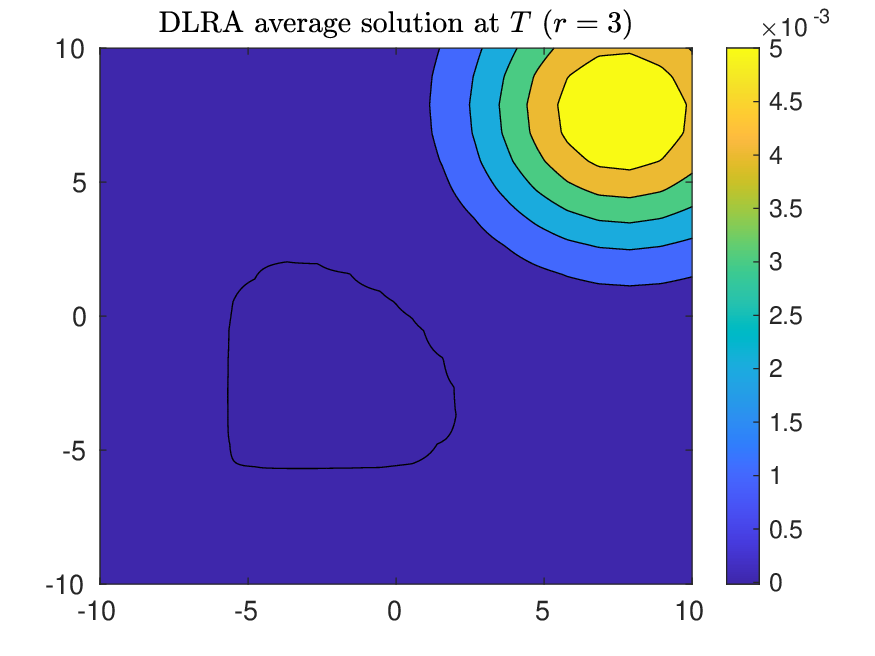}	
	\includegraphics[scale=0.33]{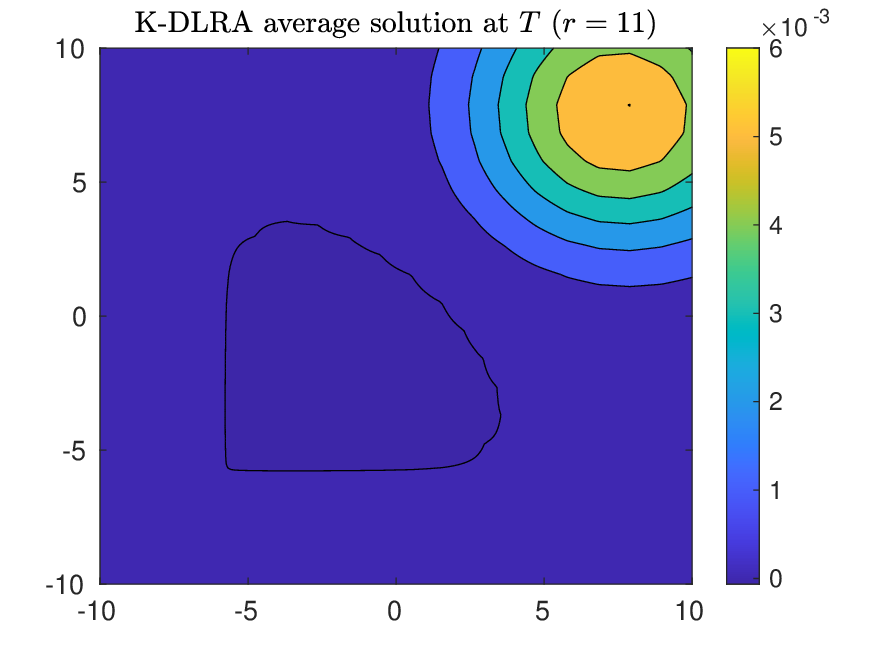}	
    \includegraphics[scale=0.33]{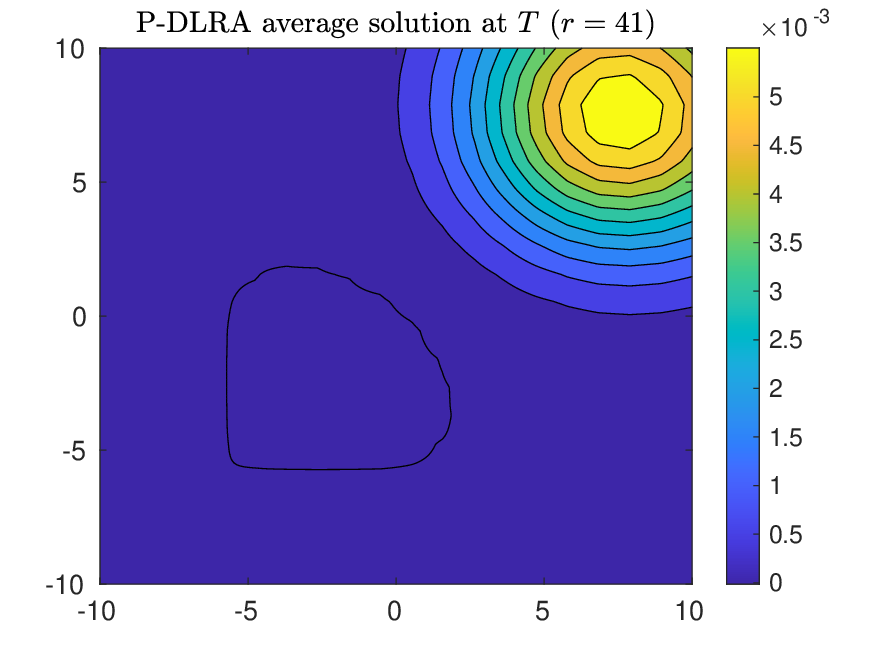}\\	

    \includegraphics[scale=0.33]{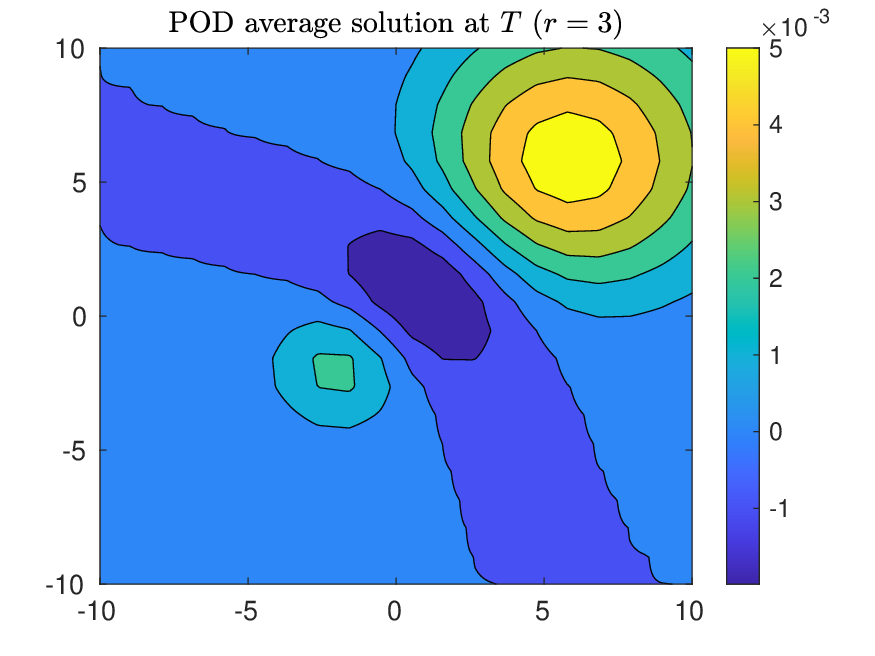}	
	\includegraphics[scale=0.33]{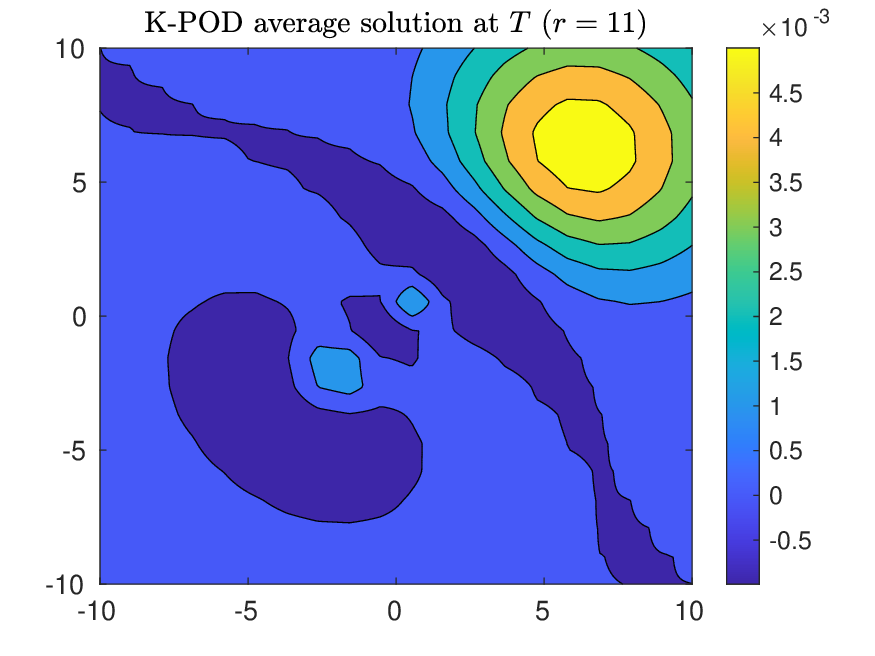}	
    \includegraphics[scale=0.33]{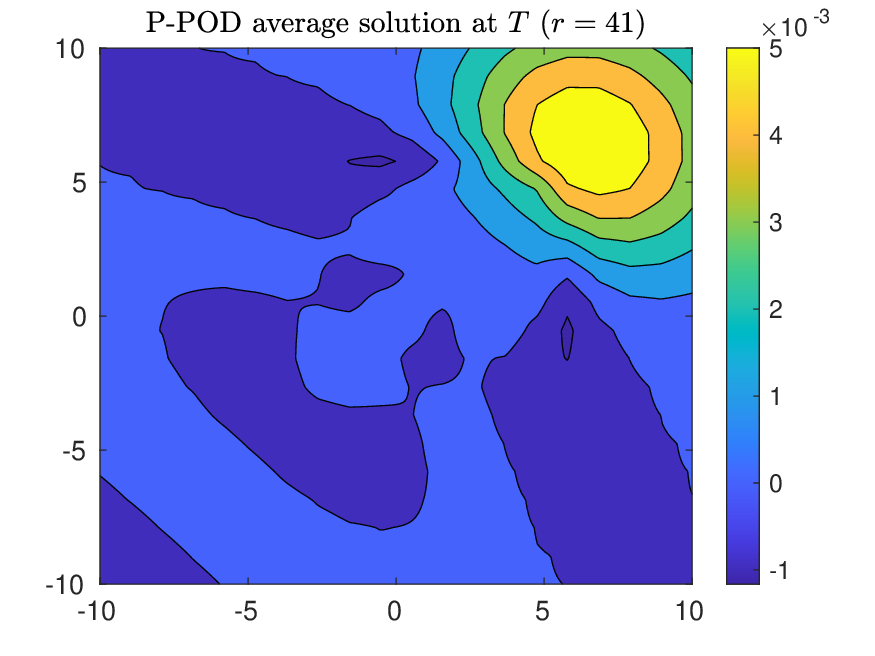}\\

\caption{Experiment 3. Solutions averaged in $\mathcal P_h$ at the final time $T$. \emph{Top}: FOM.  \emph{Middle row}: DLRA ($r=3$), K-DLRA ($r=11$), and P-DLRA ($r=41$). \emph{Bottom row}: POD ($r=3$), K-POD ($r=11$), and P-POD ($r=41$)}.
\label{pics:sol2D}
\end{figure}

\begin{figure}[H]	
\centering
	\includegraphics[scale=0.5]{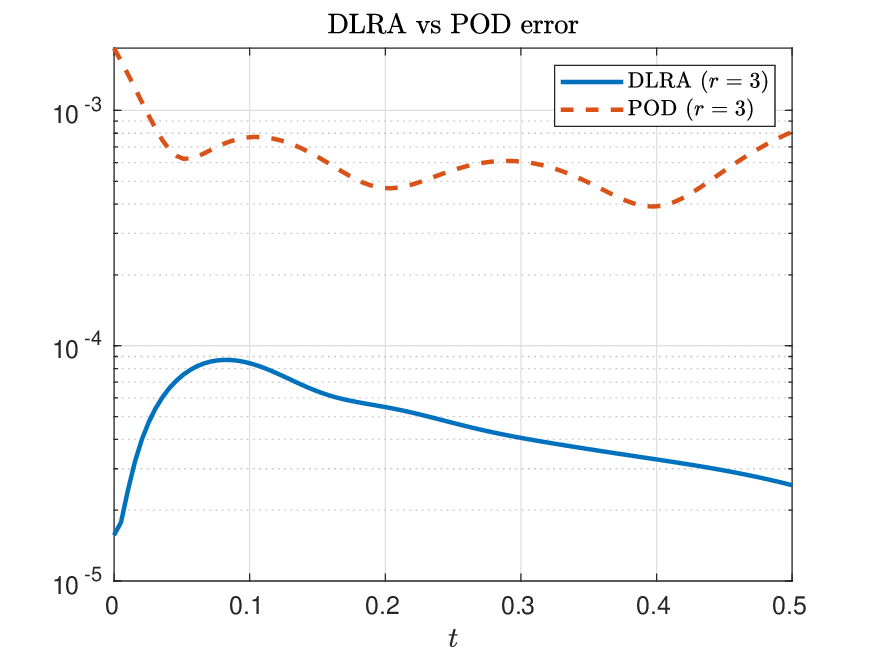}	
	\includegraphics[scale=0.5]{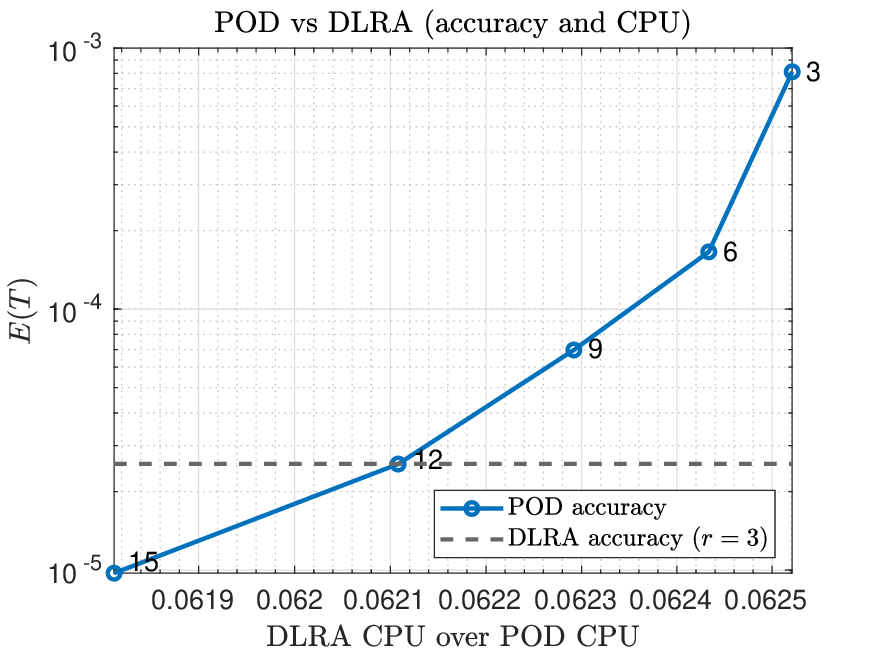}
\caption{Experiment 3. \emph{Left}: Evolution of $E(t)$ for DLRA and POD for $r=3$. \emph{Right}: POD final accuracy plot with CPU times for increasing values of $r$ to be compared to the final accuracy of the DLRA approach.}
\label{pics:DLRA_POD_2D}
\end{figure}

\begin{figure}[H]	
\centering
\includegraphics[scale=0.5]{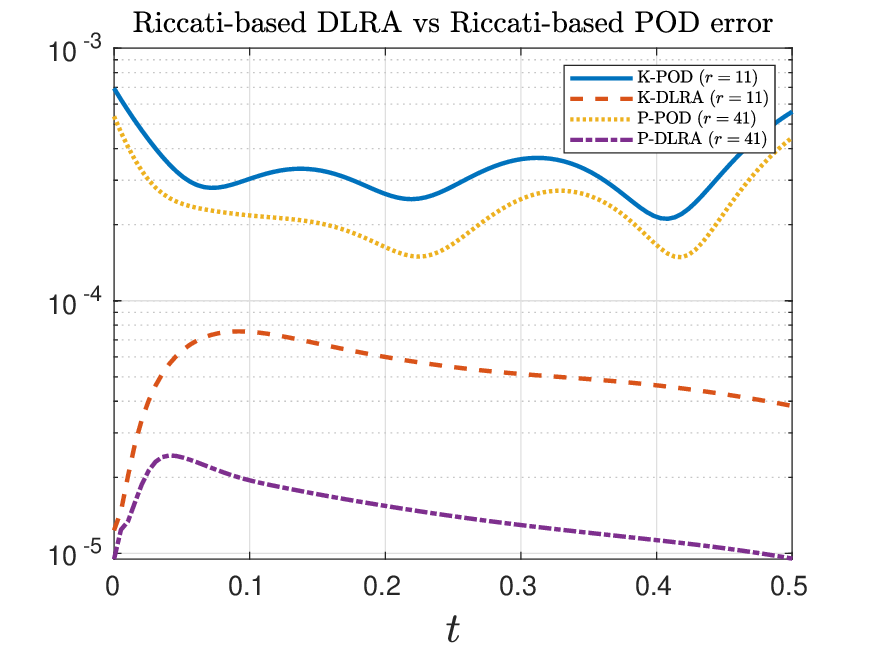}\\
	\subfloat[\centering \first{Comparison varying $r_2$}]{
    \includegraphics[scale=0.5]{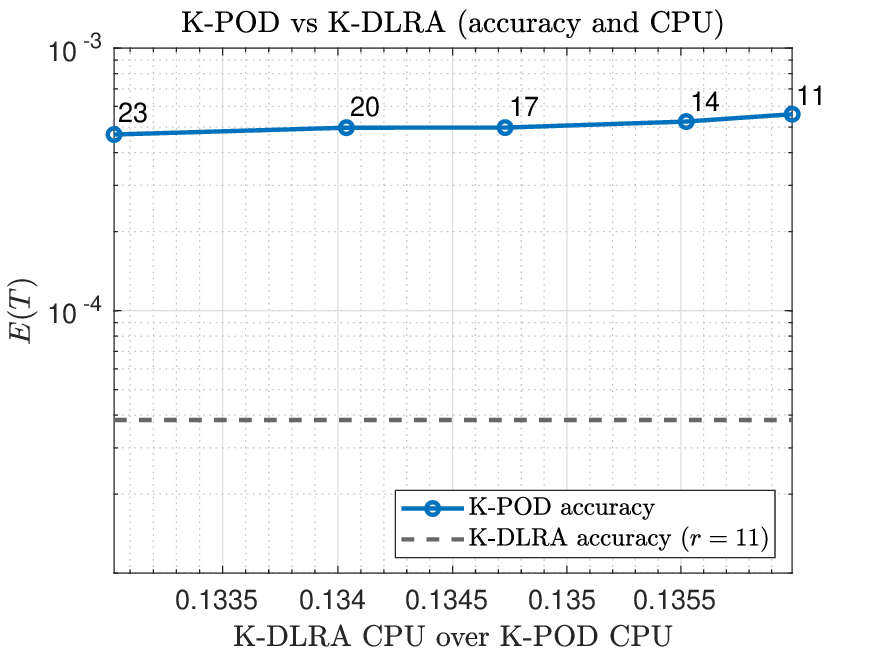}
    \label{a1b}}		
    \subfloat[\centering \first{Comparison varying $r_1$}]{
    \includegraphics[scale=0.5]{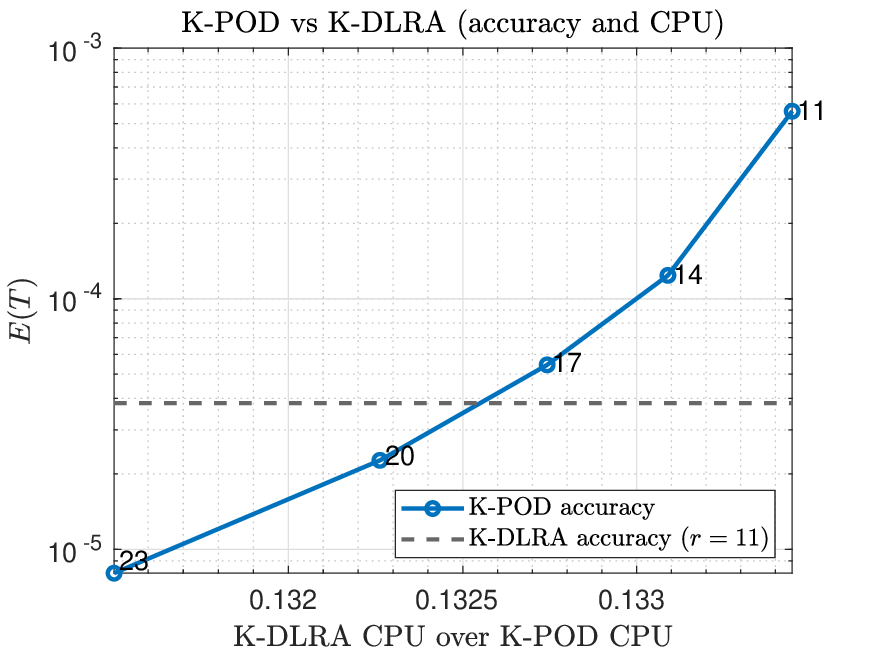}
    \label{b1b}}	\\
    \subfloat[\centering \first{Comparison varying $r_2$}]{
	\includegraphics[scale=0.5]{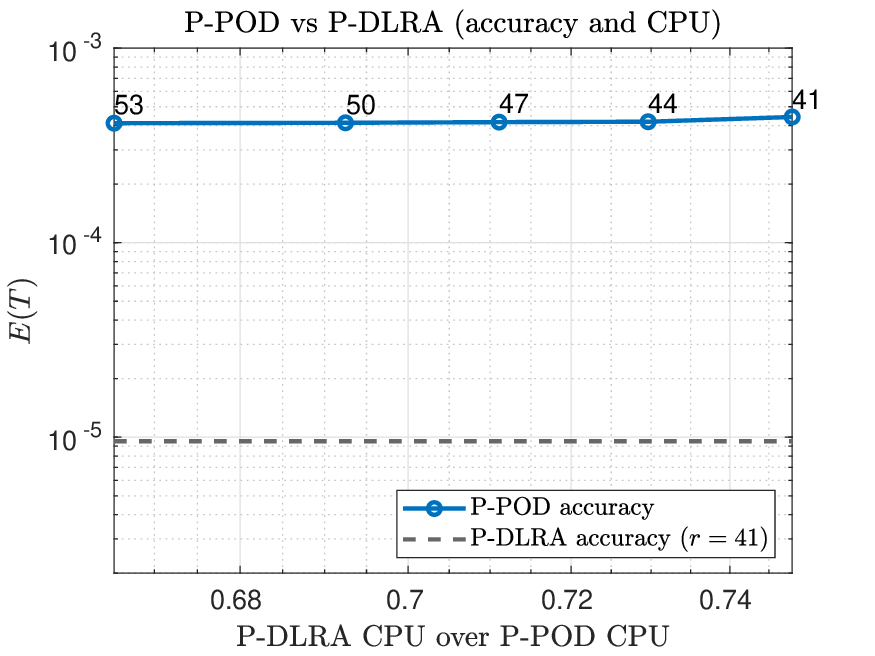}
    \label{c1b}}	
    \subfloat[\centering \first{Comparison varying $r_1$}]{\includegraphics[scale=0.5]{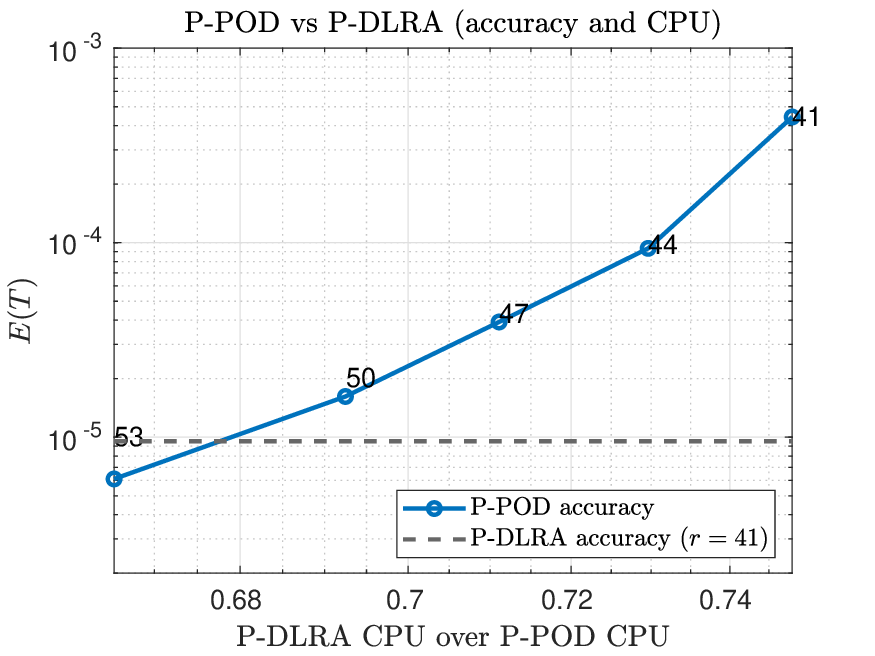} \label{d1b}}
	
\caption{Experiment 3. \emph{Top}: Evolution of $E(t)$ for K-POD, K-DLRA, P-POD, and P-DLRA with $r=11$ and $r=41$, respectively. \emph{Middle left}: K-POD final accuracy with CPU times for increasing values of $r_2$ to be compared to the final accuracy of the K-DLRA \protect\subref{a1b}. \emph{Bottom left}: analogous representation for P-POD \protect\subref{c1b}. \emph{Middle right}: K-POD final accuracy with CPU times for increasing values of $r_1$ to be compared to the final accuracy of the K-DLRA \protect\subref{b1b}. \emph{Bottom right}: analogous representation for P-POD \protect\subref{d1b}. }
\label{pics:KPerrs2D}
\end{figure}

We start with a qualitative analysis in Figure \ref{pics:sol2D}, where we plot the averaged solution in $\mathcal P_h$ at the final time. As in the one-dimensional case, DLRA and K-DLRA present an agreement with the FOM solution, but the accuracy is highly increased when the P-DLRA is used. All the POD strategies fail to represent the FOM solution, with many spurious oscillations. 
From a quantitative viewpoint, in the left plot of Figure \ref{pics:DLRA_POD_2D}, we show the evolution of $E$ in time for DLRA and POD. As expected, DLRA outperforms POD. \first{In} the right plot, we fix the DLRA final accuracy $E(T)$ for $r=3$ and we plot the final accuracy of the POD as the number of basis functions increases. To achieve comparable accuracy to DLRA, POD requires $r=12$ basis functions. Nonetheless, DLRA is unparalleled in terms of computational efficiency, being nearly two orders of magnitude faster. The CPU times and the corresponding speedups relative to the FOM solution are reported in the first column of Table \ref{tab:global_tab_2D}, where it is evident that the DLRA approach significantly outperforms the POD approach in terms of speed.
\begin{table}[H]
\caption{\centering Experiment 3. CPU comparison for DLRA, POD, K-DLRA, K-POD, P-DLRA and P-POD.}
\label{tab:global_tab_2D}
\centering
\begin{tabular}{c|cc|cc|cc}
            & \multicolumn{2}{c|}{$r=3$} & \multicolumn{2}{c|}{$r=11$} & \multicolumn{2}{c}{$r=41$} \\
            & DLRA         & POD         & K-DLRA       & K-POD       & P-DLRA       & P-POD       \\ \hline
global CPU  & $169$          & $2.71\cdot 10^{3}$         &   $370$        & $2.78\cdot 10^{3}$         &    $3.31\cdot 10^{3}$       & $4.42\cdot 10^{3}$         \\ \hline
speedup     & $374$           & $23$          &     $170$       &    $22$         &     $19$     &    $14$         \\ \hline
\end{tabular}
\end{table}
The same analysis is performed in Figure \ref{pics:KPerrs2D} for the R-DLRA and R-POD approaches. 
Firstly, the top plot demonstrates that K-DLRA offers a slight improvement over DLRA, as observed when compared to Figure \ref{pics:DLRA_POD_2D}. The middle plots present a performance analysis of the final accuracy and CPU times for K-POD, compared with the accuracy and CPU times of K-DLRA for $r=11$. Similarly to the previous test case, larger values of $r_2$ do not increase the K-POD accuracy. The only effective way to improve the accuracy is by increasing $r_1$. Notably, K-POD requires approximately twenty basis functions to be competitive with K-DLRA. However, K-DLRA remains nearly ten times faster in terms of CPU time. In the bottom plots, a similar analysis is conducted for P-POD. Increasing $r_2$ proves ineffective, while fixing $r_2=38$ and setting  $r_1=15$ lets the P-POD achieve an accuracy comparable to P-DLRA, which still outperforms the POD counterpart in terms of computational speed.

\second{

\begin{remark}
  
  We note that the maximal dimension considered in our numerical experiments is $N_h = 400$, which allows us to employ, even in the offline stage, the direct solvers available through the MATLAB functions \texttt{lyap} and \texttt{icare}. Naturally, for genuinely high-dimensional problems, these solvers become impractical, and alternative approaches must be adopted. One viable solution is offered by the MATLAB Toolbox M-M.E.S.S.~\cite{SaaKB21-mmess-3.1,benner2021matrix}, which implements low-rank solvers for large-scale symmetric matrix equations. In particular, the function \texttt{mess\_lrnm} carries out an inexact Kleinman-Newton iteration with line search for solving AREs, and allows the inclusion of an initial stabilizing feedback matrix $K(x) = R^{-1} B^\top P(x)$, which is beneficial for the cascade strategy introduced in this work. 

For instance, for the control of the two-dimensional generalized Burgers' equation, fixing the FOM dimension to $N_h = 2500$, the \texttt{mess\_lrnm} function computes an approximate solution to the ARE for a single parameter in $33.61$ seconds, achieving a Frobenius-norm residual of $6.85 \times 10^{-11}$, whereas the Newton-Kleinman method based on the \texttt{lyap} function requires $204.71$ seconds and yields a residual of $1.11 \times 10^{-12}$. Thus, in this case, the M-M.E.S.S. Toolbox achieves a speed-up of approximately sixfold, with only one order of magnitude difference in the residual. This performance advantage is expected to grow with increasing problem dimensionality.

\end{remark}

}

\begin{figure}[H]	
\centering
\includegraphics[scale=0.5]{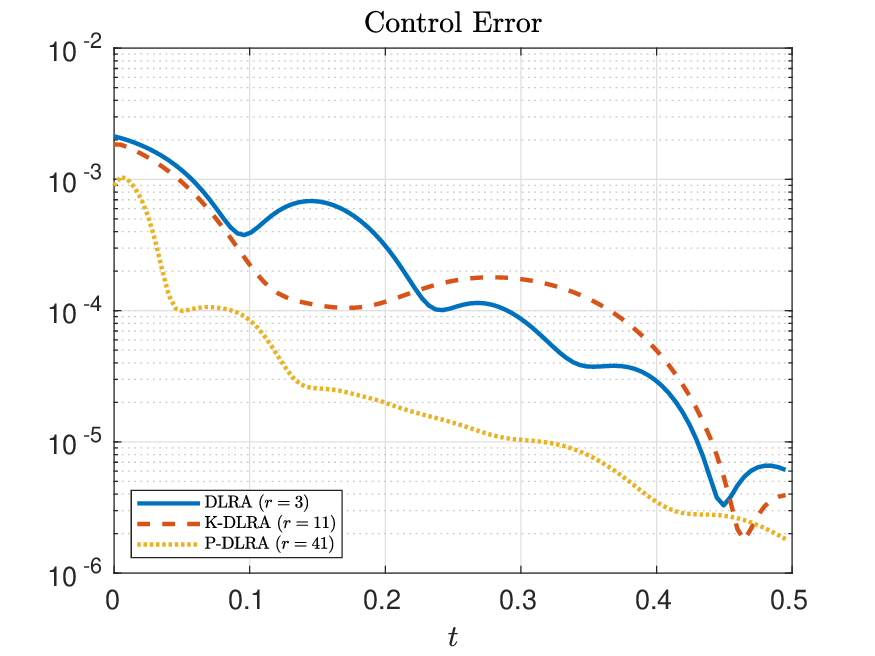}
	\includegraphics[scale=0.5]{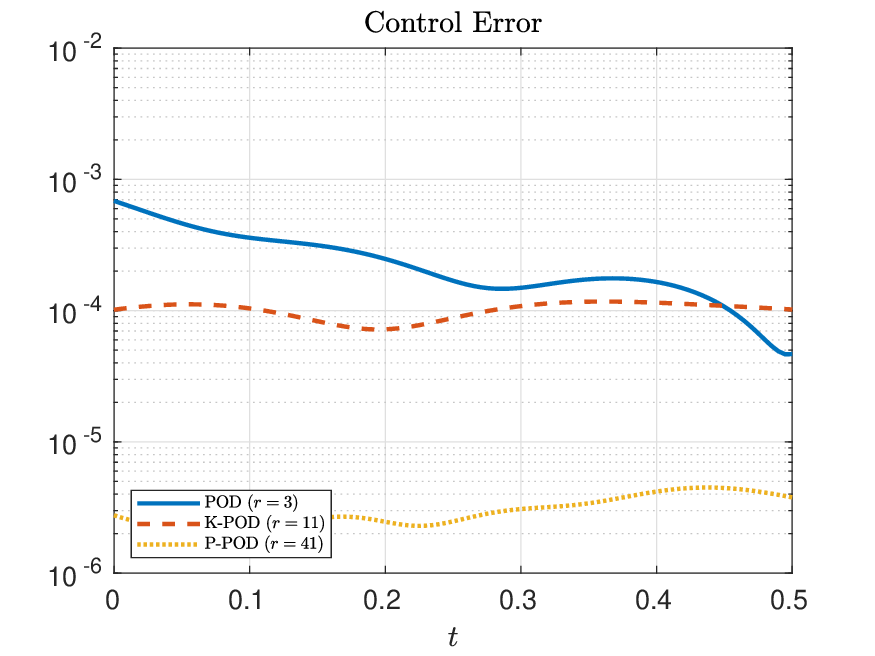}	
	
\caption{Experiment 3. Evolution of $E_c(t)$. \emph{Left plot}: DLRA, K-DRLA, and P-DLRA. \emph{Right}: analogous plot for POD, K-POD, and P-POD.}
\label{pics:control2D}
\end{figure}

Table \ref{tab:global_tab_2D} summarizes all the CPU times for all the strategies considered. 
For completeness, Figure \ref{pics:control2D} illustrates the evolution of $E_c(t)$.
On the left, we depict it for DLRA, K-DLRA, and P-DLRA. The same plot is proposed on the right for POD, K-POD, and P-POD. As in the one-dimensional case, DLRA and K-DLRA are comparable to POD and K-POD at the beginning of the simulation. The DLRA and the K-DLRA control errors drastically decrease after $t=0.3$ and start to outperform the POD and the K-POD. However, for this test case, P-POD is very accurate in terms of control error. By the end of the simulation, P-DLRA reaches the same order of accuracy. 
It is important to note that while P-POD is accurate for control errors, it remains less effective in recovering the dynamics, as shown in Figure \ref{pics:sol2D}.
Concluding, the DLRA-based approaches outperform the corresponding POD versions for all the metrics taken into consideration for the state dynamics. The P-DLRA strategy is the most accurate approach; however, as expected, it requires a greater computational time due to the large number of basis functions needed for the simulation.

\section{Conclusions}
\label{sec:conclusions}
This contribution introduces DLRA-based approaches to address parametric nonlinear feedback control under transport phenomena for stabilizing dynamical systems within the infinite horizon optimal control framework. To the best of the authors' knowledge, this is the first study to explore this direction. The DLRA framework has been enhanced through two key strategies: (i) the \texttt{C-NK} method, which accelerates the convergence of the iterative solution of the SDREs, and (ii) the leveraging of the Riccati information, leading to the development of two novel algorithms, K-DLRA and P-DLRA, that improve the accuracy of standard DLRA in control problems. These algorithms have been benchmarked against their POD counterparts in terms of accuracy with respect to the FOM solution and computational efficiency. The DLRA approaches have demonstrated superior performance in accelerating the solution of controlled dynamical systems than the FOM and the POD.

The results have been validated through one- and two-dimensional test cases involving infinite horizon feedback control for Burgers' equations. Notably, the \texttt{C-NK} strategy has been theoretically substantiated via a rigorous sensitivity analysis of the initial guesses for the \texttt{NK} iterative solver. 

This research opens several promising investigations for future work. Numerical results indicate that DLRA has challenges in recovering the FOM control, as highlighted in Remark \ref{rem:control}.
Potential future development is to combine the strengths of POD-based methods to recover control with the high precision of DLRA-based methods in the state dynamics, with the final goal of proposing an efficient algorithm with increased accuracy in representing control actions. 
Another promising improvement involves introducing adaptivity in the number of basis functions over time. Moreover, we believe that DLRA and R-DLRA might be of indisputable utility in more complicated applications in real-life scenarios, most of all in transient and convection-dominated settings.
\section*{Acknowledgments}
We acknowledge the INdAM - GNCS Project ``Metodi di riduzione di modello ed approssimazioni di rango basso per problemi alto-dimensionali" (CUP\_E53C23001670001).\\
MS thanks the ``20227K44ME - Full and Reduced order modelling of coupled systems: focus on non-matching methods and automatic learning (FaReX)" project – funded by European Union – Next Generation EU  within the PRIN 2022 program (D.D. 104 - 02/02/2022 Ministero dell’Università e della Ricerca). This manuscript reflects only the authors’ views and opinions and the Ministry cannot be considered responsible for them.
MS acknowledges the European Union's Horizon 2020 research and innovation program under the Marie Skłodowska-Curie Actions, grant agreement 872442 (ARIA).
MS thanks the ECCOMAS EYIC Grant ``CRAFT: Control and Reg-reduction in Applications
for Flow Turbulence".
The authors express their sincere gratitude for the many insightful discussions
with Stefano Massei and Cecilia Pagliantini.

\bibliographystyle{abbrvurl}
\bibliography{common_bib,biblio} 

\appendix

\section{Additional results on the sensitivity estimates}
\label{sec:appendix}
In this Appendix, we collect all the results needed to prove and validate the estimate of Section \ref{sec:sensitivity}. Section \ref{sec:ar} details the propositions necessary to the proofs of Theorems \ref{thm:NKconvergence} and \ref{thm:deltaA}. In Section \ref{sec:example}, we propose a numerical example validating the theoretical sensitivity analysis.
\subsection{Auxiliary results}
\label{sec:ar}
\begin{proposition}
    Given a vector field  $f(t, y)$ satisfying the one-sided Lipschitz condition \eqref{eq:1side_lip}, then the following estimate holds
\begin{equation}
    \Vert y(t;y_0) - y(t;\tilde{y}_0) \Vert \le \Vert y_0 - \tilde{y}_0 \Vert e^{L t}.
\end{equation}
\label{prop:decay_1side}
\end{proposition}

\begin{proof}
    Let us denote by \( y_1(t) \) and \( y_2(t) \) 
 the solutions of the ODE with corresponding initial condition $y_0$ and $\tilde{y}_0$ and
define \( \Delta y(t) = y_1(t) - y_2(t) \). Then we have:

\[ \Delta y(t) \cdot \frac{d}{dt} \Delta y(t) = \Delta y(t) \cdot (f(t, y_1(t)) - f(t, y_2(t))). \]

Using the one-sided Lipschitz condition:

\[ \Delta y(t) \cdot (f(t, y_1(t)) - f(t, y_2(t))) \leq L \| \Delta y(t) \|^2. \]

Since \( \frac{d}{dt} \| \Delta y(t) \|^2 = 2 \Delta y(t) \cdot \frac{d}{dt} \Delta y(t) \), we have:

\[ \frac{d}{dt} \| \Delta y(t) \|^2 \leq 2L \| \Delta y(t) \|^2. \]

By Gronwall's inequality, we obtain the result.
\end{proof}

\begin{proposition}
    \second{Given Assumption \ref{assum:review} \ref{hp:Acl_lip}-\ref{hp:detec}}, the vector field $f(t,y) = A_{cl}(y) y$ satisfies the one-sided Lipschitz condition \eqref{eq:1side_lip} in a domain $\Omega$ with constant $$L=\max_{x \in \Omega } \both{\textup{logm}}(A_{cl}(x)) + L_{cl} \max_{y \in \Omega} \Vert y \Vert.$$
    \label{prop:Acl_1side}
\end{proposition}

\begin{proof}
We aim to verify if the function $A_{cl}(y)y$ satisfies condition \eqref{eq:1side_lip}. Let us manipulate the left-hand side of  condition \eqref{eq:1side_lip} to obtain
$$
(x-y)^\top (A_{cl}(x)x-A_{cl}(y) y) = (x-y)^\top A_{cl}(x)(x-y)+ (x-y)^\top (A_{cl}(x)-A_{cl}(y))y \both{.}
$$
By the definition of logarithmic norm we obtain
$$
(x-y)^\top A_{cl}(x)(x-y) \le \both{\textup{logm}}(A_{cl}(x)) \|x-y\|^2\both{,}
$$
while by \second{Assumption \ref{assum:review} \ref{hp:Acl_lip}}
$$
|(x-y)^\top (A_{cl}(x)-A_{cl}(y))y| \le L_{cl} \Vert y \Vert \|x-y\|^2\both{.} 
$$
Considering the maximum in $x$ and $y$ over the domain $\Omega$ gives the result.
\end{proof}

\begin{cor}
\second{Given Assumption \ref{assum:review} \ref{hp:Acl_lip}-\ref{hp:y0_lip}}, the following estimate holds 
$$
    \Vert y(t;y_0(\mu)) - y(t;y_0(\tilde{\mu})) \Vert \le L_{y_0} \Vert \mu - \tilde{\mu} \Vert e^{Lt},
    $$
    with $L=\max_{x \in \Omega } \both{\textup{logm}}(A_{cl}(x)) + L_{cl} \max_{y \in \Omega} \Vert y \Vert$. 
    \label{cor:decay_lip_1}
\end{cor} 
\begin{proof}
By Proposition \ref{prop:Acl_1side}, we know that the function $f(t,y) = A_{cl}(y) y$ satisfies the one-sided Lipschitz condition \eqref{eq:1side_lip}. Applying Proposition \ref{prop:decay_1side} and Assumption \ref{assum:review} \ref{hp:y0_lip} we obtain the result.
\end{proof}

\begin{proposition}
\label{prop:norm_Ay}
\first{Given Assumption \ref{assum:review} \ref{hp:nabla}},
then $\Vert A_{cl}(y(t)) y(t) \Vert $ is decreasing in time.
\end{proposition}
\begin{proof}

Denoting $z(t) =A_{cl}(y(t))y(t) $, we have
$$
    \frac{d}{dt} \Vert z(t) \Vert^2 = 2 z(t))^\top \left[ (\nabla A_{cl}(y(t))) y(t) + A(y(t)) \right] z(t) \le
    $$
    $$
   2 (\both{\textup{logm}}(A(y))+\both{\textup{logm}}((\nabla A_{cl}(y)) y)) \Vert z(t) \Vert^2\both{,}
    $$
    and the result comes from Assumption \ref{assum:review} \ref{hp:nabla}.
\end{proof}

\subsection{A Numerical validation}
\label{sec:example}
In this section, we want to show how the theoretical findings derived in Section \ref{sec:sensitivity} can be applied in practice. 
Let us consider the optimal control problem driven by the following nonlinear reaction-diffusion equation:

\begin{equation}
    \frac{d}{dt} y(t,x) = \sigma \partial_{xx} y(t,x) + y(t,x)^2+ u(t,x), \quad (t,x) \in   \mathbb{R}^+ \times [0,2],
\end{equation}
        with associated cost functional
        $$
        J(u,y_0) = \int_0^{+\infty} \int_{\Omega} |y(t,x)|^2 \; dx + \int_{\Omega} |u(t,x)|^2 \; dx \, dt .
        $$
Fixing $\sigma = 10^{-3}$ and performing a semidiscretization via finite difference with $N_h=50$ grid points, we obtain
$$
\dot{y}(t) = A(y)y(t)+u(t),
$$
where
$$
A(y) = \sigma A_0 + \both{\textup{diag}}(y \circ y),
$$
with $A_0$ arising from the discretization of the Laplacian with Neumann boundary conditions, $\circ$ stands for the Hadamard product and $\both{\textup{diag}}(v)$ indicates a diagonal matrix with the components of the vector $v$ on the main diagonal.

Let us consider a set of parameterized initial conditions
$$
 y_0(x;\mu) = \mu \both{\textup{sin}}(\pi x), \quad \mu \in [1,10].
$$
We note that in this case $A(y)$ is symmetric for all $y \in \mathbb{R}^d$ and the term $F$ is an identity matrix multiplied by a scalar depending on the space discretization. As a result, the closed-loop matrix is symmetric, and all of its eigenvalues are real.
Let us denote by $\lambda_{\mu_1,\mu_2,t_1,t_2}$ the maximum eigenvalue of the matrix $A_{\mu_2,t_2}-F P_{\mu_1,t_1}$, i.e., the starting closed-loop matrix for the parameter $\mu_2$ at time $t_2$ using as initial guess $P_{\mu_1,t_1}$. Our aim is to verify whenever this value is negative and, hence, $P_{\mu_1,t_1}$ represents a stabilizing initial guess.
We proceed applying the following steps:
\begin{enumerate}
    \item we start at initial time from $\mu_1=1$ and we solve the following SDRE
$$
A^\top_{\mu_1,0}  P_{\mu_1,0} + P_{\mu_1,0} A_{\mu_1,0}-P_{\mu_1,0}F P_{\mu_1,0}  = -Q,
$$
using the Matlab function \texttt{icare};
\item we compute the approximate solution $H_{\mu_1,0}$ of the Lyapunov equation \eqref{eq:lyap_pert} by the Matlab function \texttt{lyap};
\item we note that 
$$
\Vert \Delta A_{\mu_1,\mu_2,0,0} \Vert^2 = |\mu_1^2-\mu_2^2| \Vert \both{\textup{diag}}(\both{\textup{sin}}(\pi x)) \Vert^2,
$$
where $\both{\textup{sin}}(\pi x) = [\both{\textup{sin}}(\pi x_1), \ldots , \both{\textup{sin}}(\pi x_d)]^\top$;
\item condition \eqref{condi_pert} now reads
$$
|\mu_1^2-\mu_2^2| \second{<} c = \frac{1}{4 \Vert H_{\mu_1,0} \Vert^2 \Vert \both{\textup{diag}}( \both{\textup{sin}}(\pi x)) \Vert^2},
$$
where the right-hand side is equal to $1.0011$ in this case.
Hence, condition \eqref{condi_pert} is satisfied if 
$$
\Delta \mu \second{<} - 1 + \sqrt{1+c} \approx 0.4146.
$$
\end{enumerate}

\begin{figure}[H]	
\centering
	\includegraphics[scale=0.5]{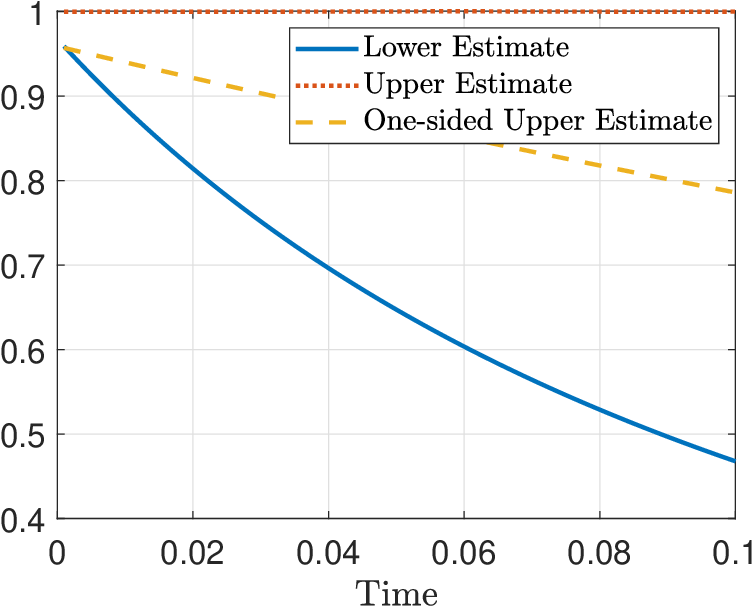}	
	\includegraphics[scale=0.5]{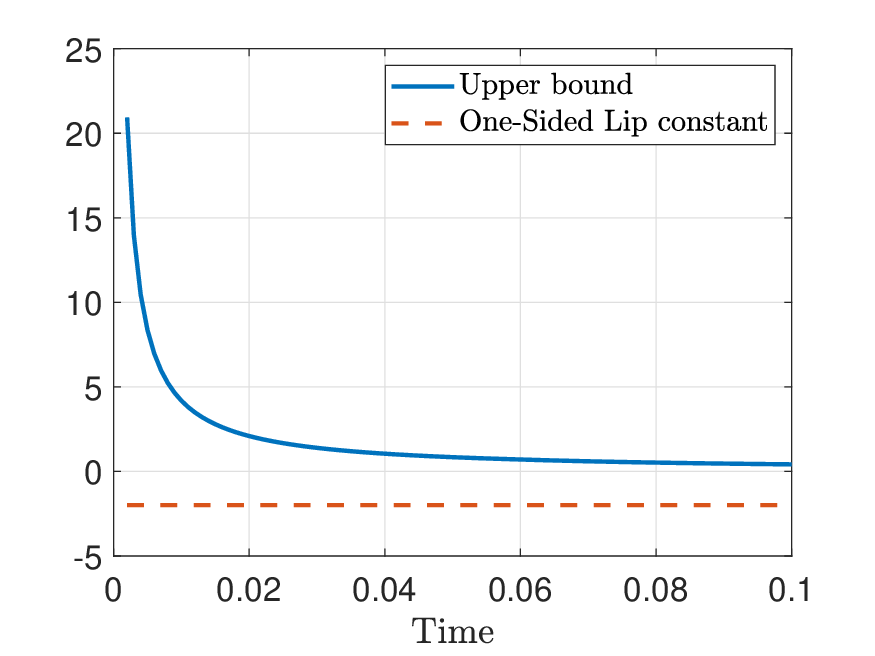}	
 \includegraphics[scale=0.5]{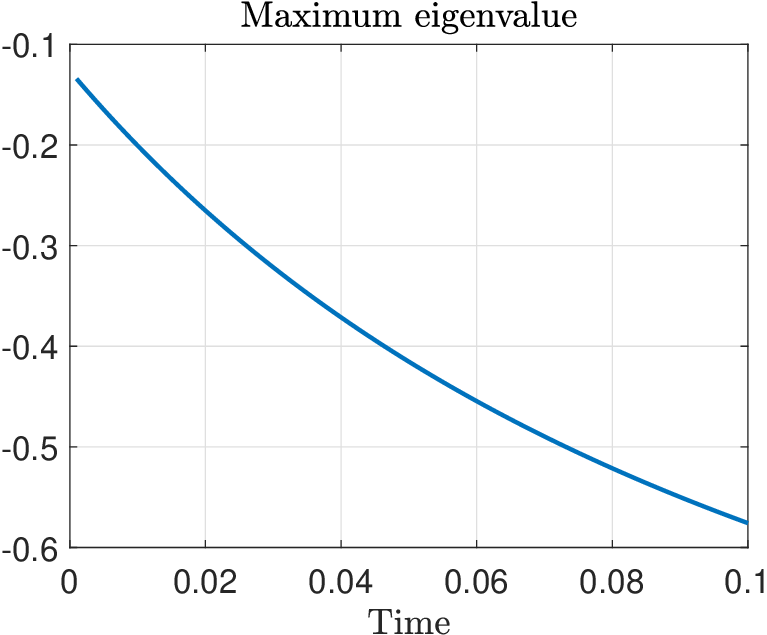}	
\caption{\emph{Top}: Verification of the estimate \eqref{condi_pert} (left) and the estimate \eqref{cond_L}  (right). \emph{Bottom}: Maximum eigenvalue $\lambda_{\mu_1,\mu_2,t,t}$ in the time interval $[0,0.1]$.}
\label{pics:estimates}
\end{figure}
If we consider $\Delta \mu = 0.45$, the condition \eqref{condi_pert} is not satisfied and we note that the maximum eigenvalue of the closed loop matrix $A_{\mu_2,0}-F P_{\mu_1,0}$ is positive ($\lambda_{\mu_1,\mu_2,0,0}=0.0019$), hence $ P_{\mu_1,0}$ is not a stabilizing initial guess. On the other hand, considering $\Delta \mu = 0.4$, as expected, $ P_{\mu_1,0}$ is a stabilizing initial guess ($\lambda_{\mu_1,\mu_2,0,0}=-0.1342
$). This provides a strategy to compute the parameter stepsize to ensure the cascade approach at the initial time.

Now let us consider the time integration with final time $T=0.1$ with $\Delta \mu = 0.1$. In the top-left panel of Figure \ref{pics:estimates}, we compare the behavior in time of the terms $g_1(t) = \Vert \Delta A_{\mu,\tilde{\mu},t,t} \Vert$ (denoted in the plot as Lower Estimate), $g_2(t) = L_A  L_{y_0} \Vert \mu - \tilde{\mu} \Vert e^{Lt}$ (One-sided Upper Estimate) and $ g_3(t) = \frac{1}{2 \Vert H_{\mu,t} \Vert}$ (Upper Estimate). The constants $L$ and $L_A$ can be estimated along the trajectories in the following way
$$
L_A = \max_{i} \frac{\Vert A(y_{\mu_1, t_i})-A(y_{\mu_2,t_i}) \Vert}{\Vert y_{\mu_1, t_i} -y_{\mu_2,t_i} \Vert},
$$
$$
L =\max_{i} \frac{(A_{cl}(y_{\mu_1, t_i})-A_{cl}(y_{\mu_2,t_i}))^\top (y_{\mu_1, t_i}-y_{\mu_2,t_i})}{\Vert y_{\mu_1, t_i} -y_{\mu_2,t_i} \Vert}.
$$
We note that the upper estimate exhibits an almost constant behavior, remaining close to values around 1.
	On the other hand, the Lower Estimate starts below as expected, since the initial choice for $\Delta \mu$ satisfies the condition \eqref{condi_pert}. The estimated one-sided Lipschitz constant $L$ is negative and equal to $-1.9897$, reflecting the decreasing behavior of the One-sided Upper Estimate and of the Lower Estimate. Condition \eqref{cond_L} is depicted in the top-right panel of Figure \ref{pics:estimates}. The right-hand side of the estimate \eqref{cond_L} is decreasing, but always positive, while the one-sided Lipschitz constant $L$ is negative. Hence, condition \eqref{cond_L} is satisfied for the entire time interval $[0,0.1]$. From this, we can infer that $P_{\mu_1,t}$ consistently serves as a stabilizing initial guess. This is demonstrated in the bottom panel of Figure \ref{pics:estimates}, where the maximum eigenvalue remains negative throughout the time interval.
 Let us now perform a similar analysis for the time step $\Delta t$. For simplicity, we will use the explicit Euler method for integration. Given the initial condition $y_{\mu_1,0} = \both{\textup{sin}}(\pi x)$,  consider the update
 $y_{\mu_1,\Delta t} = y_{\mu_1,0} + \Delta t \left ( A(y_{\mu_1,0}) - F P_{\mu_1,0}    \right ) y_{\mu_1,0}$, 
 We define $f_{cl,0}$ as $f_{cl,0} = \left ( A(y_{\mu_1,0}) - F P_{\mu_1,0}    \right ) y_{\mu_1,0}$.
 Now, we compute:
 $$
 \Vert \Delta A_{\mu_1,\mu_1,0, \Delta t} \Vert^2 = \Vert \both{\textup{diag}}(y_{\mu_1, \Delta t} \circ y_{\mu_1, \Delta t}) - \both{\textup{diag}}(y_{\mu_1, 0} \circ y_{\mu_1, 0}) \Vert = (\Delta t)^2 f_{cl,0}^2+ 2 \Delta t f_{cl,0} \cdot y_{\mu_1, 0},
 $$
 then condition \eqref{condi_pert} becomes

 $$
 (\Delta t)^2 \Vert f_{cl,0} \Vert ^2+ 2 \Delta t f_{cl,0} \cdot y_{\mu_1, 0} - \frac{1}{2 \Vert H_{\mu_1,0} \Vert^2}  \second{<} 0,
 $$
which leads to the solution
 $$
\Delta t \second{<}\frac{ - f_{cl,0} \cdot y_{\mu_1, 0} + \sqrt{(f_{cl,0} \cdot y_{\mu_1, 0})^2+\frac{1}{2 \Vert H_{\mu_1,0} \Vert^2 }}}{\Vert f_{cl,0} \Vert^2}  \approx  1.5694.
$$
In fact, for $\Delta t = 1.6$, the maximum eigenvalue of the closed loop matrix $A(y_{\mu_1, \Delta t})-F P_{\mu_1,0}$ is positive ($\lambda_{\mu_1,\mu_1,0,\Delta t}=  0.8752$), whereas for $\Delta t = 1.5$, it is negative ($\lambda_{\mu_1,\mu_1,0,\Delta t}=  -0.7858$).
This section justifies the use of the cascade information in solving the SDRE by means of \texttt{NK}.

\end{document}